\documentclass[12pt,a4paper]{article}
\usepackage{amsmath, amssymb, subfigure}
 \usepackage{graphicx}
\usepackage{hyperref}

\newtheorem{theorem}{Theorem}[section]
\newtheorem{lemma}[theorem]{Lemma}
\newtheorem{proposition}[theorem]{Proposition}
\newtheorem{definition}[theorem]{Definition}

\newtheorem{remark}[theorem]{Remark}

\allowdisplaybreaks

\def\neweq#1{\begin{equation}\label{#1}}
\def\endeq{\end{equation}}

\def\Om {{\Omega}}
\def\la {{\lambda}}

\def \Ab{{\bf A}}
 
\def\Rb {{\bf R}} 

\def \Inte{{\rm Int\,}}

\newcommand {\pa}{\partial}

\newcommand {\beq}{\begin{equation}}
\newcommand{\eeq}{\end{equation}}
\numberwithin{equation}{section}
\newenvironment{proof}{\begin{trivlist}
                       \item[]\hspace{0cm}{\bf Proof: }
                       \hspace{0cm} }{\hfill $\Box$
                     \end{trivlist}}

\begin{document}

\title{\bf\Large On spectral  minimal partitions II,\\  the case of the rectangle}
\author{V. Bonnaillie-No\"el\footnote{IRMAR, ENS Cachan Betagne, Univ. Rennes I, CNRS}, 
B. Helffer\footnote{D\'epartement de Math\'ematiques, Universit\'e Paris-Sud}, 
T. Hoffmann-Ostenhof\footnote{Institut f\"ur Theoretische Chemie, Universit\"at Wien
{\it and} International Erwin Schr\"odinger Institute for Mathematical Physics}}
\date{\today}

\maketitle

\begin{abstract}
In continuation of \cite{HHOT}, we discuss the question of 
spectral minimal $3$-partitions for the rectangle $]-\frac a2,\frac a2[\times ] -\frac b2,\frac b2[\,$, 
with $0< a\leq b$. It has been observed in \cite{HHOT} that when 
$0<\frac ab < \sqrt{\frac 38}$ the minimal $3$-partition is obtained by the three 
nodal domains of the third eigenfunction corresponding to the  three 
rectangles  $]-\frac a2,\frac a2[\times ] -\frac b2,-\frac b6[$,  
$]-\frac a2,\frac a2[\times ] -\frac b6,\frac b6[$ and 
$]-\frac a2,\frac a2[\times ] \frac b6, \frac b2[$.  We will describe a 
possible mechanism of transition for increasing  $\frac ab$ 
between these  nodal minimal $3$-partitions and non nodal minimal 
$3$-partitions at the value $ \sqrt{\frac 38}$ and discuss the existence 
of symmetric candidates for giving minimal $3$-partitions when $ \sqrt{\frac 38}<\frac ab
  \leq 1$.  Numerical analysis leads very naturally to nice questions
 of isospectrality which are solved by introducing Aharonov-Bohm Hamiltonians 
 or by going on the double covering of the punctured rectangle.\\[5pt]
\end{abstract}

\section{Introduction}
In continuation of \cite{HHOT}, we have analyzed in \cite{HH:2006}
 the question of minimal $3$-partitions for the disk
 and introduced new tools for this partially successful analysis.  
In the same spirit,  we discuss here  the similar question  for the rectangle $\Rb_{a,b}:=]-\frac a2,\frac
  a2[\times ] -\frac b2,\frac b2[$, with $0< a\leq b$. For a given
 partition\footnote{See the next section for precise definitions.}  $\mathcal D$  of an open set $\Omega$ by $k$ open subsets $D_i$, we can consider
\begin{equation}
\Lambda(\mathcal D)=\max_{i=1,\ldots,k}\la(D_i)\,,
\end{equation}
where $\la(D_i)$ is the ground state energy of the Dirichlet Laplacian on $D_i\,$. 
We denote the infimum on every $k$-partitions of $\Omega$ by
\begin{equation}
\mathfrak L_k(\Omega)=\inf_{\mathcal D\in \mathfrak O_k}\Lambda(\mathcal D)\,.
\end{equation} 
We look for minimal $k$-partitions, that is partitions 
such that  $\mathfrak L_k(\Omega) =\Lambda(\mathcal D)\,$.\\

It has been observed in \cite{HHOT} that, when $0< \frac ab <  \sqrt{\frac 38}$ the minimal 
$3$-partition is given by the three nodal domains of the third 
eigenfunction corresponding to the  three rectangles  $]-\frac a2,\frac
  a2[\times ] -\frac b2,-\frac b6[\,$,  $]-\frac a2,\frac
  a2[\times ] -\frac b6,\frac b6[$  and  $]-\frac a2,\frac
  a2[\times ] \frac b6, \frac b2[\,$.

In the case when $ \sqrt{\frac 38}<\frac ab \leq 1\,$,  we can show
as for the disk, see \cite{BHV},  that
$\mathfrak L_3(\Rb_{a,b})$ is not an eigenvalue. Indeed,
$\la_2(\Rb_{a,b})=\la_3(\Rb_{a,b})<\la_4(\Rb_{a,b})$. By Theorem~\ref{L=L} below,
see \cite{HHOT} for the proof, $\mathfrak L_3(\Rb_{a,b})$ cannot be an
eigenvalue
 and hence the associated minimal partition cannot be nodal.\\

 We will describe in Section~\ref{sec.transition} a possible 
mechanism of transition for increasing  $\frac ab$ between these 
nodal minimal $3$-partitions and non nodal minimal $3$-partitions at the value 
$ \sqrt{\frac 38}$ and discuss the existence of symmetric candidates for 
giving minimal $3$-partitions when $ \sqrt{\frac 38}<\frac ab \leq 1\,$ in Sections~\ref{sec.isosquare}-\ref{sec.heuristics}. 

\begin{figure}[h!]
\begin{center}
\hfill\subfigure[First candidate $\mathcal D$]
	{\includegraphics[width=4cm]{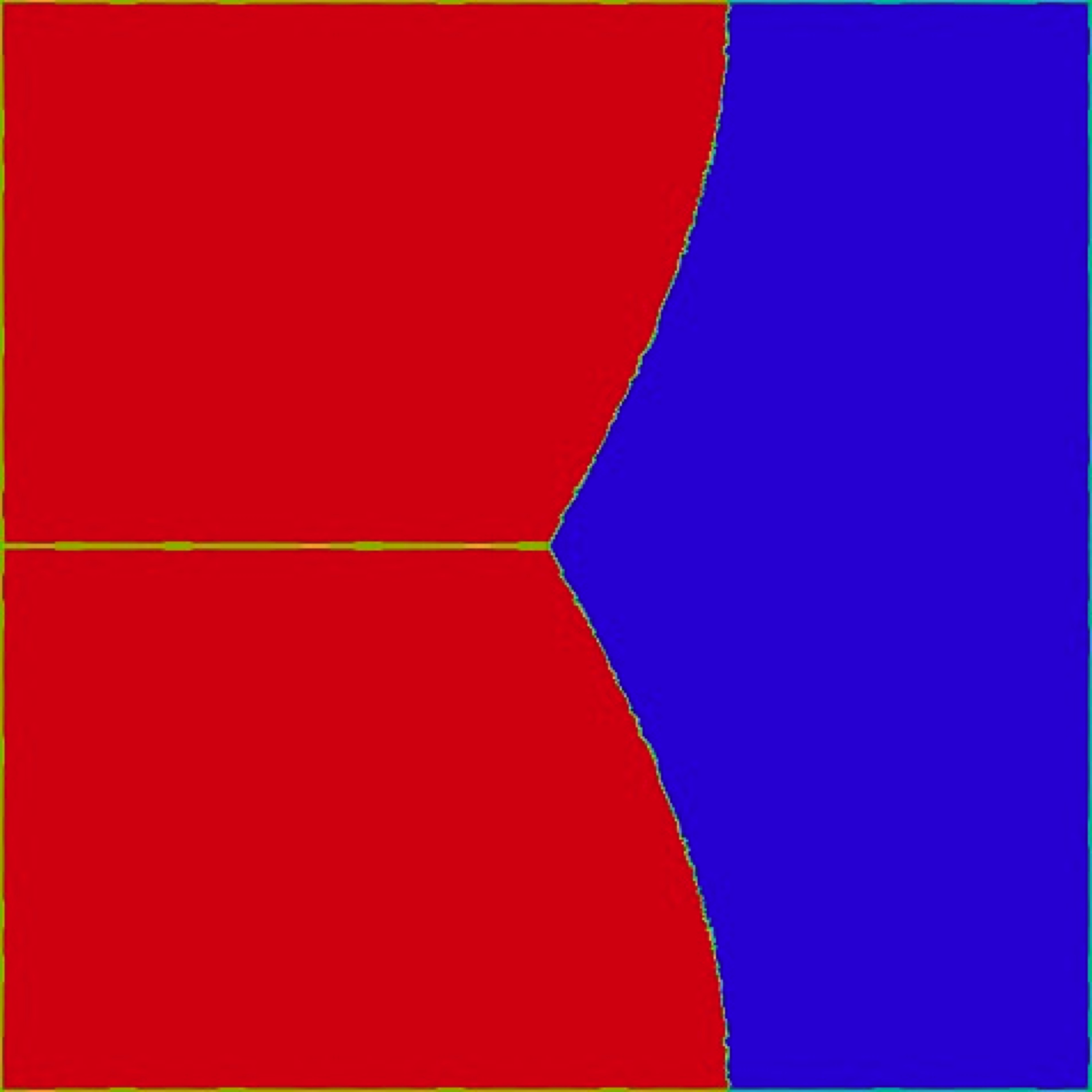}\label{picturesquare}}
\hfill\subfigure[Second candidate $\mathcal D^{new}$]
	{ \includegraphics[width=4cm]{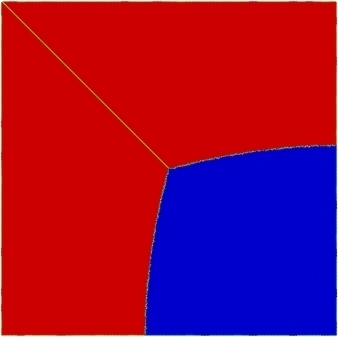}\label{diagsquare}}\hfill \ 
\caption{Candidates for the minimal  $3$-partition of the square.}
\end{center}
\end{figure}

We can exhibit numerically some candidates for the minimal $3$-partition using symmetry.
Assuming that there is a minimal partition which is symmetric with respect to the 
axis $\{y=0\}\,$, and intersecting the partition with the half-square  
$]-\frac 
12,\frac 12[\times ]0,\frac 12[\,$,  one is reduced to analyze a family of 
Dirichlet-Neumann problems. Numerical computations\footnote{See 
{\sf
  http://w3.bretagne.ens-cachan.fr/math/simulations/MinimalPartitions/} } performed by V.~Bonnaillie-No\"el and G.~Vial  (in January 2006) lead to a natural candidate  $\mathcal D$ for a symmetric minimal partition (see Figure~\ref{picturesquare}).
We observe numerically that the three lines of $N(\mathcal D)$ (i.e. the interior boundary of the subdomains $D_{i}$ of the partition, see Definition \eqref{eq.NmathcalD}) meet at the center $(0,0)$ of the 
square. As expected by the theory they meet at this critical point with equal angle $\frac{2\pi}{3}$
 and meet  the boundary orthogonally. 
This choice of symmetry is not unique. By exploring numerically 
the possibility of a minimal partition with another symmetry (diagonal), we get the surprise to find 
another candidate $\mathcal D^{new}$ with  $\Lambda(\mathcal
D^{new})=\Lambda(\mathcal D)$ (see Figure \ref{diagsquare}).\\
This  leads very naturally to nice questions of isospectrality 
which are solved using Aharonov-Bohm Hamiltonians or by 
going on the double covering of the punctured rectangle. Sections 
\ref{sec.isospec}-\ref{sec.isosquare} concern these questions. 
 
The paper is organized as follows: \\
In Section~\ref{sec.2}, we recall some notation and properties concerning the nodal partitions. \\
In Section~\ref{srectangle}, we start with the  analysis  of the
rectangle based on results of \cite{HHOT} and enumerate in particular
all the possible Courant-sharp situations. There appears the limiting
case $a/b=\sqrt{3/8}$ detailed in Section~\ref{sec.transition} in
which 
 a mechanism is proposed which explains the  transition between the nodal $3$-partition and the non nodal one.\\
Motivated by the numerical simulations for the square
we find two candidates for the minimal 3-partition with the same energy.
We analyze in Section~\ref{sec.isospec} the Aharonov-Bohm Hamiltonian and give some isospectral properties for rectangles. This theory is applied in Section~\ref{sec.isosquare} to explain numerical simulations for the minimal $3$-partitions on the square.
The paper ends in Section~\ref{sec.heuristics} with some heuristics on the deformation of symmetric minimal partitions corroborated with numerical simulations for rectangles $]-\frac{\epsilon\pi} 2,\frac{\epsilon\pi} 2[\times ]-\frac\pi 2,\frac\pi2[$ from $\epsilon=\sqrt{\frac 38}$ to $\epsilon=1$.
 
\section{Definitions, notation and previous results\label{sec.2}}
As in \cite{HHOT} (see also \cite{HH:2006}),
 we consider the Dirichlet Laplacian 
on a bounded domain $\Om\subset \mathbb R^2$, which is piecewise
 $C^\infty$. We denote, for any open domain $D$,  the lowest 
eigenvalue of the Dirichlet realization $H(D)$ of  $-\Delta $ in $D$
 by $\la(D)$.
For any function $u\in C_0^0(\overline\Om)\,$, we introduce 
\begin{equation}
N(u)=\overline{\{x\in \Om\:\big|\: u(x)=0\}}
\end{equation}
and call the components of $\Om\setminus N(u)$ the nodal domains of $u$. The number of 
nodal domains of such a function is denoted by $\mu(u)$. 

We now recall the main definitions and results concerning spectral
minimal partitions and refer to \cite{HHOT} for proofs and
details. For $k\geq 1$ ($k\in \mathbb N$), we call a $k$-partition of $\Om$  a family $\mathcal D=\{D_i\}_{i=1}^k$ of pairwise  disjoint 
open  domains in $\Omega$. 
Such a  partition  is called  \textbf{strong} if 
\begin{equation}\label{regcov}
\Inte(\overline{\cup_{i=1}^k D_i})\setminus \pa \Om=\Om\,.
\end{equation}
We denote by $\mathfrak O_k$ the set of such partitions.\\
For $\mathcal D\in \mathfrak O_k$ we introduce 
\begin{equation}\label{LA}
\Lambda(\mathcal D)=\max_{i=1,\ldots,k}\la(D_i)\,,
\end{equation}
and 
\begin{equation}\label{Lfrak}
\mathfrak L_k=\inf_{\mathcal D\in \mathfrak O_k}\Lambda(\mathcal D)\,.
\end{equation} 
$\mathcal D $ is called 
 a (spectral)\footnote{We will omit the word spectral.}  minimal $k$-partition if $\mathfrak
 L_k(\Omega)=\Lambda(\mathcal D)$
 and a nodal minimal $k$-partition if $\mathcal D$ consists of the
 nodal domains of an eigenfunction of $H(\Omega)$.\\
To each strong partition $\mathcal D$ we associate a graph $G(\mathcal D)$ in the following way:\\
We say $D_i,D_j\in \mathcal D$ are \textbf{neighbors}, and we denote this by $D_i\sim D_j$, if 
\begin{equation}\label{DsimD}
\Inte(\overline{D_i\cup D_j})\setminus \pa \Om\text{ is connected}.
\end{equation}
We associate to each $D_i\in \mathcal D$ a vertex $v_i$ and for each pair $D_i\sim D_j$ 
an edge $e_{i,j}$. This defines  a planar graph $G(\mathcal D)$. We
say that the partition is {\bf admissible} if the corresponding graph is bipartite.
 We recall that a nodal partition is always admissible.\\ 
Attached to a  partition $\mathcal D$,  we can associate a closed set $N\in \overline\Om$
defined by 
\begin{equation}\label{eq.NmathcalD}
N(\mathcal D)=\overline{\bigcup_i(\pa D_i\cap \Om)}\,.
\end{equation}
This leads us to introduce the  set $\mathcal M(\Om)$ of regular closed
sets $N$ which share with nodal sets all the standard properties except at isolated critical points where they have  only\footnote{We do not assume anymore that the number of lines
 arriving at a critical point is even.} the \textbf{equal angle meeting property}. 
   More precisely, we recall from \cite{HHOT} the following definition.
\begin{definition}~\\
A closed set $N\subset \overline\Omega$ belongs to $\mathcal M(\Omega)$ if $N$ satisfies
\begin{enumerate}
\item There are finitely many distinct $x_{i}\in \Omega\cap N$ and associated positive integers $\nu(x_{i})\geq 3$ such that, in a sufficiently small neighborhood of each $x_{i}$, $N$ is the union of $\nu(x_{i})$ smooth arcs (non self-crossing) with one end at $x_{i}$ and such that in the complement of these points in $\Omega$, $N$ is locally diffeomorphic to a regular curve. The set of these critical points of $N$ is denoted by $X(N)$.
\item $\partial\Omega\cap N$ consists of a finite set of points $z_{i}$ such that at each $z_{i}$, $\rho(z_{i})$ arcs hit the boundary with $\rho(z_{i})\geq 1$. 
We denote the set of critical points of $N\cap\partial\Omega$  by $Y(N)$.
\end{enumerate}
\end{definition}
A partition $\mathcal D$ 
is called regular if the corresponding $N(\mathcal D)$ is regular. Let
us now recall the main theorems.

\begin{theorem}\label{thstrreg}~\\
For any $k$ there exists a minimal regular strong $k$-partition
 and any minimal  $k$-partition admits a representative which is 
regular and strong.
\end{theorem}

The existence of a minimal regular strong partition has
been shown\footnote{See also \cite{He} and \cite{CL}.}  by Conti-Terracini-Verzini in \cite{CTV0, CTV2,
  CTV:2005}, while the second part of the theorem has
been shown in \cite{HHOT}.\\
In the following, we always consider the regular representative without
mentioning it explicitly.
We have now the following converse theorem, see \cite{HHOT}.
\begin{theorem}\label{partnod}~\\
Assume that there is a minimal admissible  $k$-partition. Then this  partition is associated to the nodal set of an 
eigenfunction corresponding to $\mathfrak L_k(\Omega)$.
\end{theorem}
This result was completed in \cite{HHOT} in relation with the Courant-sharp property. 
We recall that if  $u$ is an eigenfunction of the Dirichlet Laplacian in $\Omega$
 attached to the  $k$-th eigenvalue $\lambda_k$, then 
 Courant's Theorem says  that the number of nodal domains $\mu(u)$ satisfies
 $\mu(u)\leq k\,.$ Pleijel's Theorem says that, when the dimension is
$\geq 2$, then the previous inequality is strict for $k$ large.\\
As in 
\cite{AHH:2004}, we say  that  $u$ is \textbf{Courant-sharp} if $\mu(u)=k$.
For any integer $k\ge 1$, we denote by $L_k(\Omega)$ the smallest
eigenvalue
for which the  eigenspace 
contains an eigenfunction with $k$ nodal domains. In general we have
\beq \label{compades3}
\lambda_k(\Omega)\leq\mathfrak L_k(\Omega)\leq L_k(\Omega)\,.
\eeq
The next result of \cite{HHOT} gives  the full picture of the equality cases:
\begin{theorem}\label{L=L}~\\
If $\mathfrak L_k(\Omega)=L_k(\Omega)$ or $\lambda_k(\Omega) =\mathfrak L_k(\Omega)$, then $$\la_k(\Omega)=\mathfrak L_k(\Omega)=L_k(\Omega)\,.$$
 In addition,  any minimal (regular) $k$-partition is a nodal
 partition corresponding  to an  eigenfunction associated  to $\la_k(\Omega)$. 
\end{theorem}

As a consequence of Euler's Formula, we have described in \cite{HH:2006}
 the possible topological types of a non admissible minimal
 $3$-partition
 of a connected regular open set.
 
\begin{figure}[h!t]
\begin{center}
\hfill \includegraphics[height=3cm]{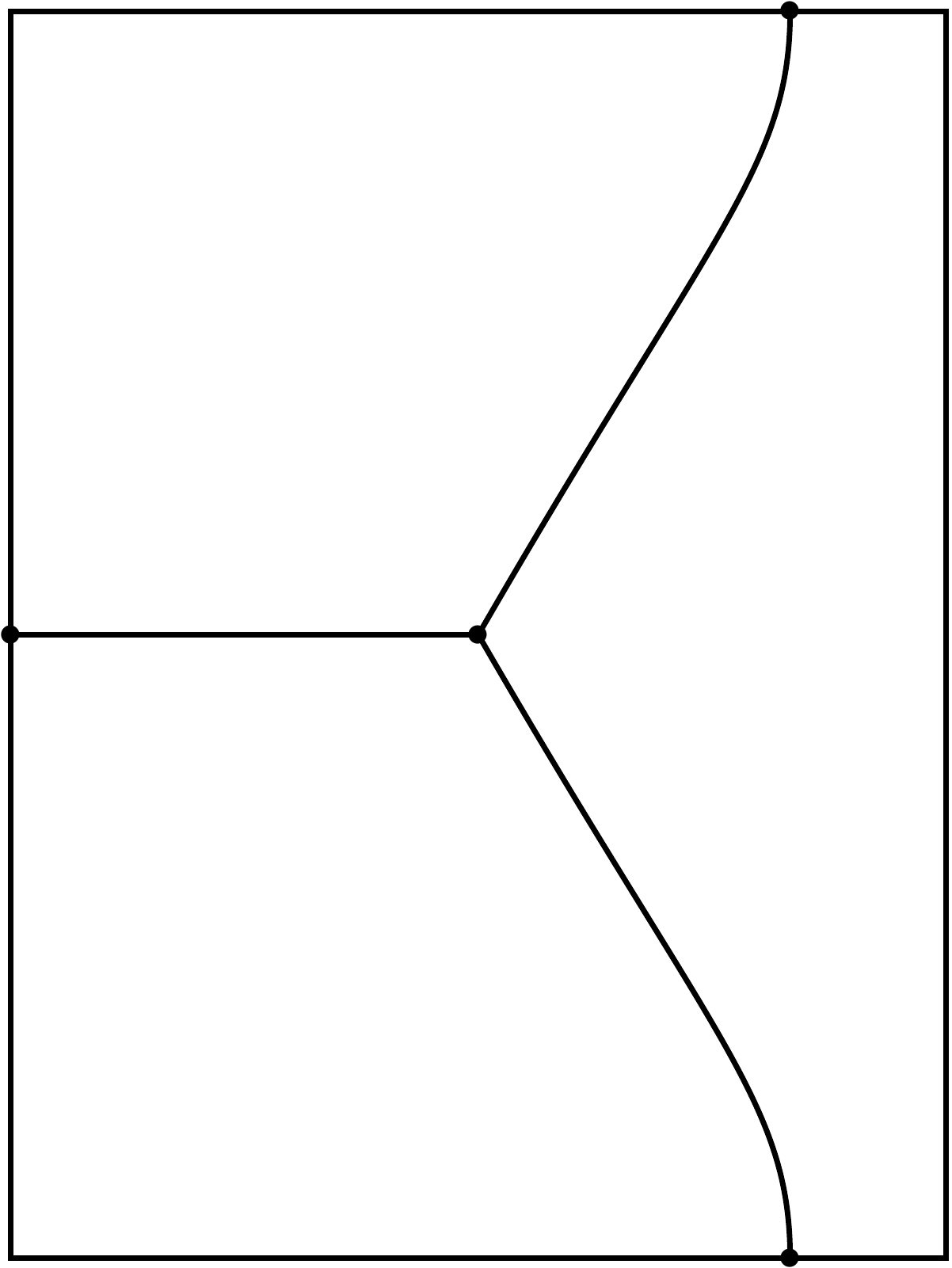} \hfill
\includegraphics[height=3cm]{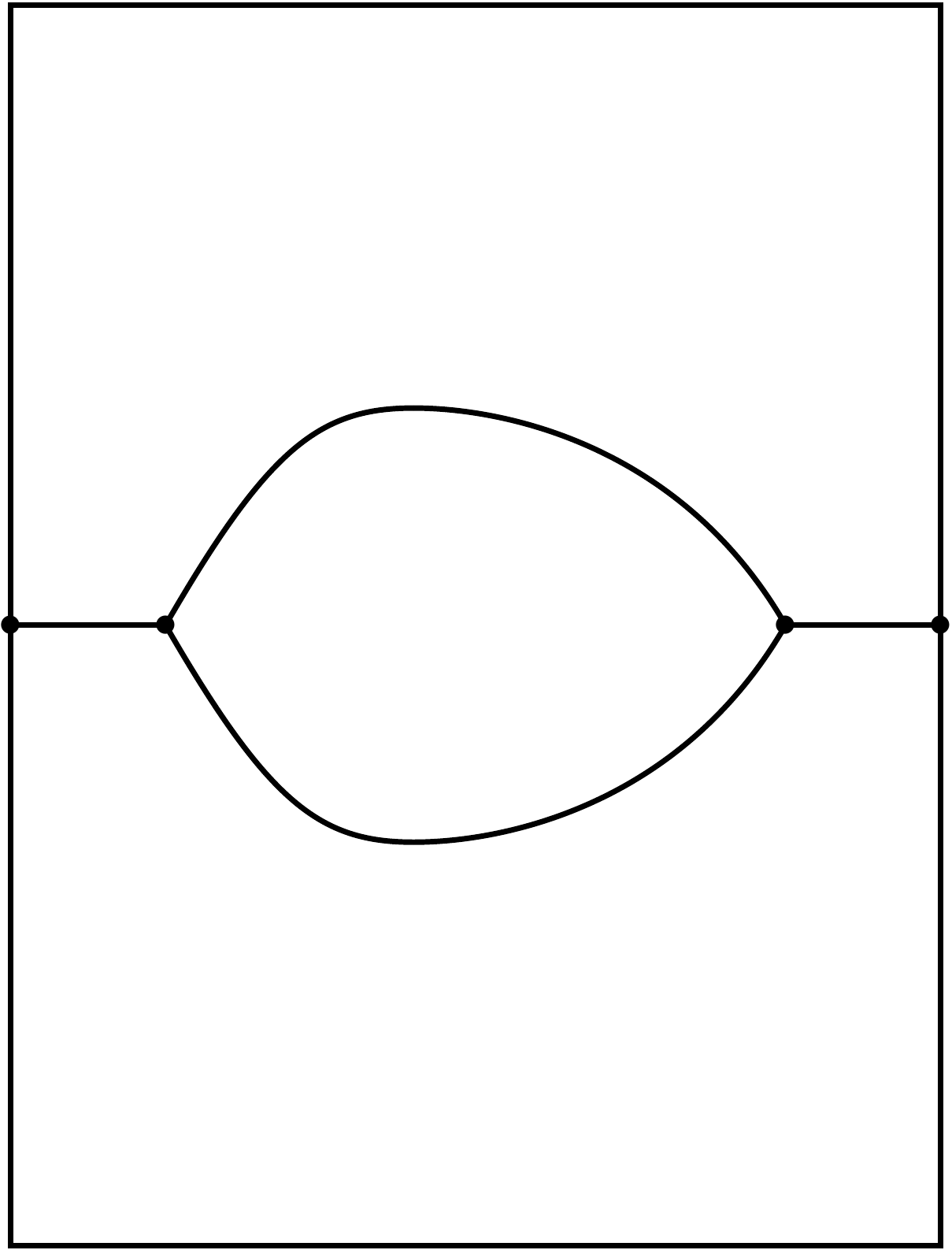} \hfill
\includegraphics[height=3cm]{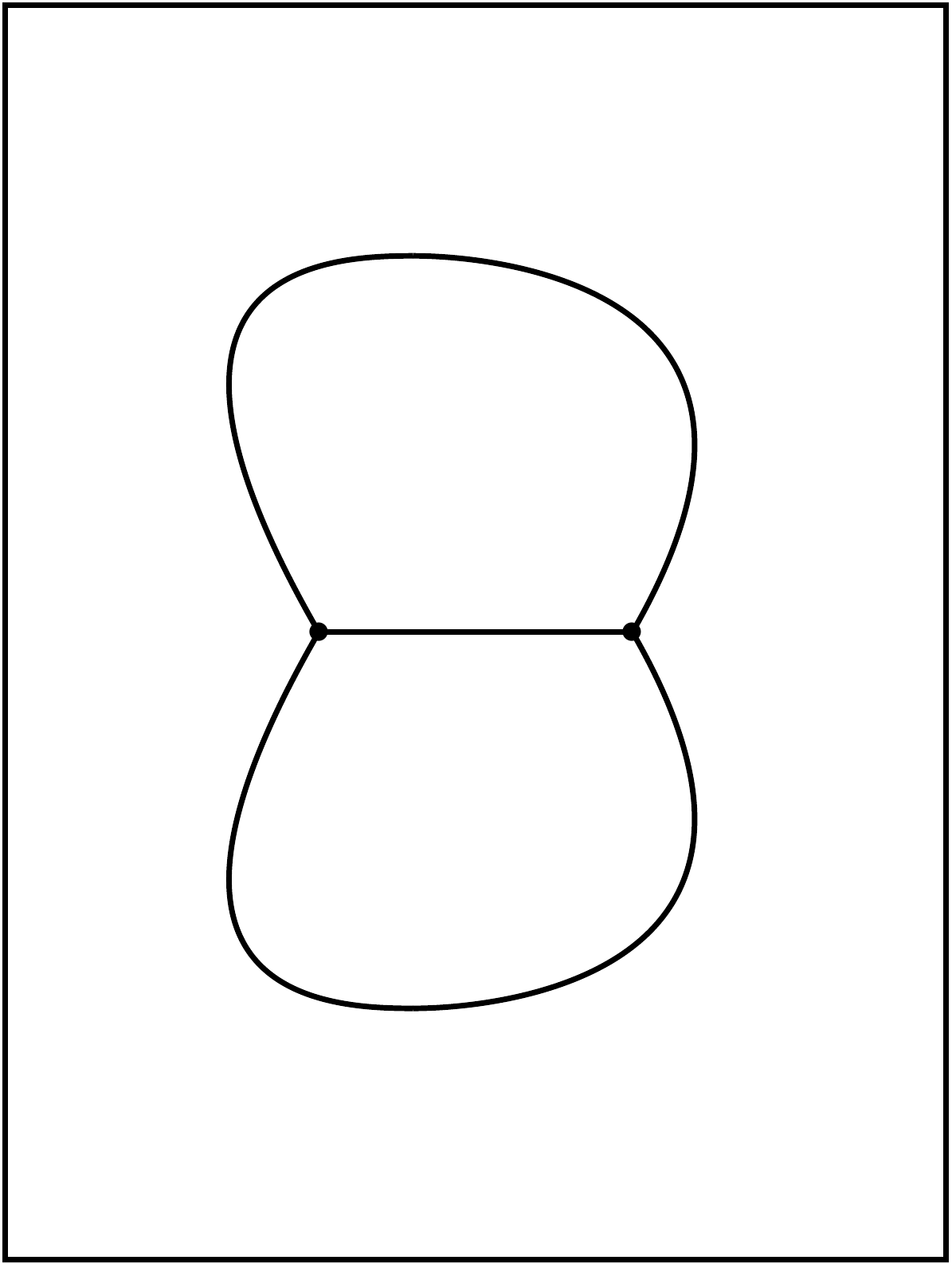} \hfill \ 
\caption{The three configurations (a), (b) and (c), with an additional
  symmetry with respect to the $x$-axis.\label{fig.typeabc}}
\end{center}
\end{figure}

\begin{proposition}\label{mp3}~\\
Let $\Omega$ be simply connected and let us consider a minimal $3$-partition $\mathcal D=(D_1,D_2,D_3)$ of
$\Omega$ associated
to $\mathfrak L_3 (\Omega)$. Let us  suppose that 
\begin{equation}\label{nbp}
\lambda_3(\Omega) <\mathfrak L_3(\Omega)\,.
\end{equation} 
For  $N =N(\mathcal D) $, we denote by $\nu(x_{i})$ and $\rho(z_{i})$ the number of arcs associated with $x_{i}\in X(N)$, respectively $z_{i}\in Y(N)$. 
Then there are three possibilities (see Figure~\ref{fig.typeabc}): \\
\textbf{(a)} $X(N)$ consists of one point $x$ with $\nu(x)=3$ and 
 $Y(N)$ consists of either 
three distinct $y_1,\,y_2,\,y_3$ points with $\rho(y_1)=\rho(y_2)=\rho(y_3)=1\,$, two distinct points $y_1,\,y_2$ with 
$\rho(y_1)=2$, $\rho(y_2)=1$ or one point $y$ with $\rho(y)=3\,$.\\
\textbf{(b)}
$X(N)$ consists of two distinct points $x_1,\,x_2$ with
$\nu(x_1)=\nu(x_2)=3\,$.\break  $Y(N)$ consists either of two points $y_1,\,y_2$
such that $\rho(y_1)+\rho(y_2)=2$ or of  one point $y$ with $\rho(y)=2\,$.\\
\textbf{(c)}
$X(N)$ consists again of two distinct points $x_1,\,x_2$ with $\nu(x_1)=\nu(x_2)=3$, but $Y(N)=\emptyset\,$. 
\end{proposition}

\section{A first  analysis  of the rectangle}\label{srectangle}
This part is taken from \cite{HHOT}. Note that when $\Omega$ is  a
rectangle $\Rb_{a,b}:= ]-\frac a2,\frac a2[\times ]-\frac b2,\frac b2[$, 
the spectrum of $H(\Rb_{a,b})$ and the properties of the eigenfunctions are analyzed
 as toy models  in \cite[\S4]{Pleijel:1956}. The spectrum is given by
$$
\lambda_{m,n}:= \pi^2 \left(\frac{m^2}{a^2} + \frac{n^2}{b^2}\right)
\quad\mbox{ with }  (m,n)\in (\mathbb N^*)^2\,.$$
These eigenvalues are simple
 if ${a^2}/{b^2}$ is irrational. 
Except for specific  remarks concerning the square or the case when $a^2/b^2=3/8$, we assume 
\begin{equation}\label{irra}
 a^2/b^2 \mbox{ is  irrational}.
\end{equation}
 So we can associate to each eigenvalue
$\lambda_{m,n}$, up to a non-zero multiplicative constant,  a unique eigenfunction $u_{m,n}$
 such that $\mu(u_{m,n}) = mn$.
Given $k\in \mathbb N^*$, the lowest eigenvalue corresponding to $k$
nodal domains is given, at least under Assumption \eqref{irra},  by
\begin{equation}
L_k(\Rb_{a,b}) = \pi^2 \inf_{mn=k} \left(\frac{m^2}{a^2} + \frac{n^2 }{b^2}\right)\,.
\end{equation}
In the case when $a^2/b^2$ is rational we could have problems in
 the case of multiplicities. We have then to 
analyze a continuous families of nodal sets of eigenfunctions living
 in an eigenspace of dimension $>1$. 
We will see  in
 Section~\ref{sec.transition} 
 that it is just for these values that new nodal  partitions
 may appear, which could be, by deformation,  
the starting point of non admissible minimal partitions.\\
We now recall all the possible Courant-sharp situations\footnote{We do not know whether for
certain $a^2/b^2$ rational additional Courant sharp eigenvalues could
show up.}
 described in \cite{HHOT}:
\begin{enumerate}
\item 
 $m=3$, $n=2$ and 
$\displaystyle
\frac 35  \leq \frac{a^2}{b^2} \leq \frac 58 \,.
$
\item  $m=2$, $n=2$ and 
$\displaystyle
\frac 35 \leq \frac{a^2}{b^2} \leq 1 \,.
$
\item  $m=1$, $n>1$ and 
$\displaystyle
\frac{a^2}{b^2} \leq  \frac{3} {n^2-1}  \,.
$
\end{enumerate}
If we now focus on the case $k=3$, we get that $\lambda_3(\Rb_{a,b})$
is Courant-sharp iff $a^2/b^2 \leq  3/8$. Hence,  the limiting situation is 
$$
\frac{a^2}{b^2} =  \frac 38\,.
$$
This corresponds to a double eigenvalue and to the pairs
$(m,n)=(1,3)$ and $(m,n)= (2,1)$.

\section{Transition from Courant-sharp to a non nodal minimal partition \label{sec.transition}}
We start from a rectangle with $a= \pi \epsilon $ and $b=\pi $ and would like
 to analyze $\mathfrak L_3 (\epsilon):=\mathfrak L_{3}(\Rb_{\pi\epsilon,\pi})$. 
 The critical situation corresponds to
\begin{equation}
\epsilon = \sqrt{ 3/8}\,.
\end{equation}
So the first result (deduced from \cite{HHOT}) which was recalled
 in the previous section writes:
\begin{proposition}\ 
\begin{enumerate}
\item If $\epsilon \leq \sqrt{ 3/8}$, 
 then 
$
\mathfrak L_3(\epsilon) = 9 + {1}/{\epsilon^2}\,
$
and $\mathfrak L_3(\epsilon)$ is an eigenvalue.
\item
If $ \sqrt{3/8}< \epsilon \leq 1$, then
$ 
\mathfrak L_3(\epsilon) < 9 + {1}/{\epsilon^2}\,.
$
\end{enumerate}
\end{proposition}
We would like to understand how the transition
 occurs at $\epsilon = \sqrt{3/8}$.
We make the assumption that in the deformation the 
minimal partition remains symmetric with respect to 
$y=0$. This is indeed the case for $\epsilon < \sqrt{3/8}\,$, because 
the eigenfunction corresponding to $\lambda_{1,3}$ is 
$\cos \frac x \epsilon\; \cos 3  y\; $ and the corresponding 
nodal lines are composed of two horizontal lines $y=- \pi/ 6$ and $y= {\pi}/{6}$. This is also 
the case for $\epsilon = \sqrt{3/8}$ because all the 
eigenfunctions have this symmetry and any minimal partition is nodal.\\

The numerical computations for the square (see Figure~\ref{picturesquare} and \cite[\S 3.2]{BHV}) 
push us to conjecture that the nodal lines $N(\mathcal D)$ for the minimal 
$3$-partition $\mathcal D$ is  for $\epsilon > \sqrt{ {3}/{8}}$
 is the union of a segment $[-\pi/2,x_0(\epsilon)]$ (on the line $y=0$)
 and of two symmetric arcs connecting the point $(x_0(\epsilon),0)$ to the
 boundary of the rectangle (up and down).\\

The first conjecture is that $x_0(\epsilon)$ is increasing monotonically 
from $-\pi/2$ to $0$ for $\epsilon \in ]\sqrt{ {3}/{8}}, 1]$. This has 
been partially verified (admitting the symmetry and that the minimal
 $3$-partition is of type (a)) numerically.\\
The second conjecture  is that the minimal $3$-partition will ``tend''
 as $\epsilon$ tends to $\sqrt{ {3}/{8}}$ from above
 to a nodal partition, losing there its non bipartite character.\\

The point here is that, when $\epsilon =  \sqrt{{3}/{8}}$, 
 we have an eigenvalue of multiplicity two  giving the possibility
 of constructing a continuous family of nodal  minimal
 $3$-partitions. For this, we consider the family 
$$
\varphi_{\alpha,\beta} (x,y) = \alpha  \cos \frac x\epsilon \cos 3y +
\beta \sin \frac{2x}{\epsilon}  \cos y \;,
$$
with $\alpha^2+ \beta^2\neq 0$, and analyze their zero set. Of course,
for $t\neq 0$, $\varphi_{\alpha,\beta}$ and $\varphi_{t\alpha,t\beta}$
 have the same zero set.\\
 
\begin{figure}[h!]
\begin{center}
\subfigure[$\alpha=1$, $\beta=0$]{\includegraphics[height=4cm]{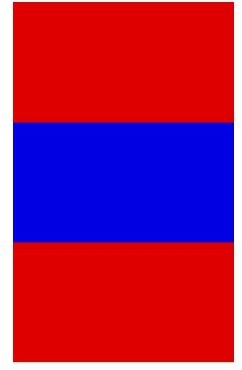}}
\subfigure[$\alpha=5$, $\beta=1$]{\includegraphics[height=4cm]{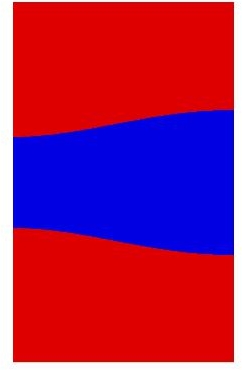}}
\subfigure[$\alpha=2$, $\beta=1$]{\includegraphics[height=4cm]{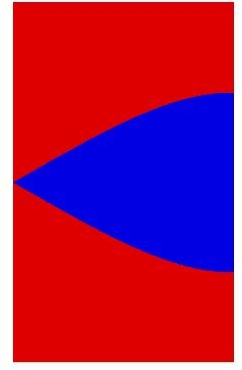}}
\subfigure[$\alpha=1$, $\beta=2$]{\includegraphics[height=4cm]{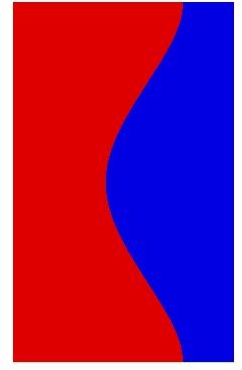}}
\subfigure[$\alpha=0$, $\beta=1$]{\includegraphics[height=4cm]{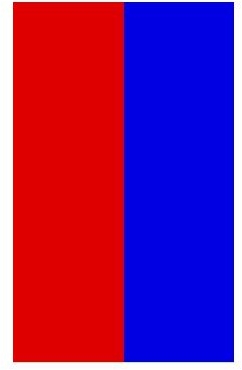}}
\caption{Nodal sets of $\varphi_{\alpha,\beta}(x,y) = \alpha \cos \frac x \epsilon \,\cos 3y\,+ \beta
\sin \frac{2x}{\epsilon}\, \cos y$.\label{fig.nodalsets}}
\end{center}
\end{figure}

We  first show that the zero set of $\varphi_{\alpha,\beta}$ has 
no critical point inside the rectangle for the critical value of $\epsilon$.
Using the factorization of  $\varphi_{\alpha,\beta}$ in the form
$$
\varphi_{\alpha,\beta} (x,y)= \cos y \cos \frac{x}{\epsilon}
 \left(\alpha  \left(1- 4 \sin^2 y \right)  + 2 \beta \sin 
\frac{x}{\epsilon}\right) \;,
$$
we observe that $\varphi_{\alpha,\beta}=0$ is equivalent inside the
rectangle to $\psi_{\alpha,\beta}=0$ with
\beq
\psi_{\alpha,\beta} (x,y) :=\alpha \left( 1- 4 \sin^2 y \right) + 2 \beta \sin \frac x \epsilon\,.
\eeq
Hence we now look at the (closure of the) zero set of
$\psi_{\alpha,\beta}$ in the rectangle and particularly 
to the critical points inside and at the  boundary. 
If we look at the zeros of  $\pa_x \psi_{\alpha,\beta} $, we obtain
$$ 
 \beta \cos  \frac{x}{\epsilon} =0\;.
$$
This implies $\beta =0$ and we get that $\varphi_{\alpha,0}(x,y) =\alpha \cos
3y\cos\frac x\epsilon\;$. 
Hence we get  that, for any $(\alpha,\beta)\neq (0,0)$, there
is no  critical point inside this  rectangle.\\
It remains to look at  what is going on at the boundary  and to determine
the singular points where two lines touch.
An analysis of the function $(x,y)\mapsto \alpha (1-4\sin^2y) + 2
\beta\sin \frac x\epsilon$ at the boundary shows that critical 
points at the boundary can only 
occur for $y=0$ and $\alpha \pm 2 \beta =0$.\\
Hence we have obtained that the
 only nodal sets having critical sets are
 (up to a multiplicative constant) the nodal domains of the
 eigenfunctions $\varphi_{2,1}$ and $\varphi_{2,-1}$.\\
Figure~\ref{fig.nodalsets} gives the nodal set of the functions $\varphi_{\alpha,\beta}$ 
for several values of $(\alpha,\beta)$.

\section{Aharonov-Bohm  Hamiltonian and isospectrality\label{sec.isospec}}
As explained in the introduction, 
this new analysis is motivated by numerical computations showing that
pushing the same idea as in \cite{BHV}  but using the symmetry
 with respect to  the diagonal, one gets the same eigenvalue
 and again a $3$-partition
 with unique singular point at the center. This will be further explained in
 more detail in Subsection \ref{sssq}. 

\subsection{Basic material}
This material appears already in \cite{HHOO} and is motivated by the
 work of Berger-Rubinstein \cite{BeRu}. 
If $\Omega$ is an open set such that $0\in \Omega$, 
a possibility  is to consider  the Aharonov-Bohm  Laplacian
in  the punctured  $\dot{\Omega}=\Omega \setminus \{0\}$,
with the singularity of the potential at the center
and normalized  flux $\Phi=1/2$. 
The magnetic potential with flux $\Phi$
 takes the form
\begin{equation}
{\bf A}(x,y) = (A_1(x,y),A_2(x,y))=\Phi\, \left( -\frac{y}{r^2}, \frac{x}{r^2}\right)\,.
\end{equation}
We know that the magnetic field vanishes identically in
  $\dot{\Omega}$ 
  and, in any cut domain (such that it becomes simply connected), 
 one has
\begin{equation}
A_1 \,dx + A_2\, dy  = \Phi \, d \theta\,,
\end{equation}
where
\begin{equation}
z= x+ i y = r \exp i \theta\,.
\end{equation}
So the Aharonov-Bohm operator in any open set $\dot \Omega \subset \mathbb
R^2\setminus \{0\}$ will always be defined  by considering the Friedrichs
extension starting from $C_0^\infty(\dot \Omega)$
 and the associated differential operator is
\begin{equation}
-\Delta_{{\bf A}} := (D_x - A_1)^2 + (D_y-A_2)^2\,.
\end{equation}
In polar coordinates (which of course are not very 
 well adapted to the square but permit
 a good analysis at the origin), the
Aharonov-Bohm Laplacian reads:
\begin{equation}
- \Delta_\Ab = (D_x + \frac 12 \, \frac{\sin \theta}{r})^2 +
 (D_y -\frac 12 \,\frac{\cos \theta}{r})^2\,,
\end{equation}
or
\begin{equation}\label{BAsecondform}
- \Delta_\Ab = -\frac{\pa^2}{\pa r^2} - \frac 1r \frac{\pa}{\pa r}
 + \frac{1}{r^2} (i\pa_\theta + \frac 12)^2\,.
\end{equation}
This operator is preserving  ``real'' functions  in some modified  sense. 
Following \cite{HHOO}, we will say that a function $u$ is  K-real, if
it  satisfies
\begin{equation}\label{defreel}
K u =u\,,
\end{equation}
where $K$ is an antilinear operator in the form
\begin{equation}\label{defK}
K = \exp i \theta \, \Gamma  \;,
\end{equation}
and where $\Gamma$ is the complex conjugation 
\begin{equation}\label{defGamma}
\Gamma u = \bar u\,.
\end{equation}
The fact that $-\Delta_\Ab$ preserves K-real eigenfunctions
 is an immediate consequence of 
\begin{equation}\label{commu}
K \circ (-\Delta_\Ab) = (-\Delta_\Ab) \circ K\,.
\end{equation}
\begin{remark}~\\
Note that our choice of $K$ 
 is not unique. $K_\alpha = \exp i \alpha K$  is also antilinear,
 satisfies \eqref{commu}  and 
$$
K_\alpha^2 = Id\,.
$$
\end{remark}
As observed in \cite{HHOO}, it is easy to find a basis of K-real eigenfunctions. 
These eigenfunctions (which can be identified with real antisymmetric eigenfunctions of the 
Laplacian on a suitable double  covering $\dot{\Omega}^{\mathcal R}$ of 
$\dot{\Omega}$) have a nice nodal structure (which is locally 
in the covering the same as the nodal set of real eigenfunctions 
of the Laplacian), with the specific property that the number of 
lines in $\dot{\Omega}$ ending at the origin should be odd. More generally 
a path of index one around the origin should always meet an odd number of nodal lines.

\subsection{Symmetries of the rectangle}
We consider now a domain $\dot \Omega$ which has the symmetries of a
rectangle. More precisely, if  we denote by $\sigma_1$
 and $\sigma_2$ the symmetries respectively defined by
\beq\label{symrect}
\sigma_1(x,y)= (-x,y)\;,\;
\sigma_2(x,y)= (x,-y)\;,
\eeq
we assume that
\beq
\sigma_1 \dot \Omega =\dot \Omega\;,\; \sigma_2 \dot \Omega =\dot
\Omega\,.
\eeq
We assume that $\Omega$ is convex (to simplify)  and write
$$
\Omega \cap \{y=0\} = ]-\frac a 2,\frac a2[\times \{0\}\;,
$$
and
$$
\Omega \cap \{x=0\} = \{0 \}\times ]-\frac b 2,\frac b2[\,.
$$
If $\Sigma_1$ is the natural action on $L^2(\dot \Omega)$ associated with
 $\sigma_1$ 
\begin{equation}\label{defSigma1}
\Sigma_1 u (x,y)=u(-x,y)\,,
\end{equation}
we observe that  the  Aharonov-Bohm  operator 
does not commute with
 $\Sigma_1$ but with the antilinear operator 
 \begin{equation}\label{defSigma1c}
 \Sigma_1^c:= i\, \Gamma \Sigma_1\,.
\end{equation}
So if $u$ is an eigenfunction, $\Sigma_1^c u$ is an
eigenfunction.\\
Moreover, and this explains the choice of ``$i$'' before $\Gamma
\Sigma_1$, since  $K$ and $\Sigma_1^c$ commute,
\begin{equation}\label{commSigma1}
K\circ \Sigma_1^c =  \Sigma_1^c \circ K\;,
\end{equation}
$ \Sigma_1^c u$
 is also a  K-real  eigenfunction   if $u$ is a K-real  eigenfunction.  
One can do the same thing with $\Sigma_2$, associated with 
 $\sigma_2$, 
\begin{equation}\label{defSigma2}
\Sigma_2 u(x,y)= u (x,-y)\,.
\end{equation}
This leads  this time to 
\begin{equation}\label{defSigma2c}
\Sigma_2^c = \Gamma \Sigma_2\,.
\end{equation} 
Similarly, we have
\begin{equation}\label{commSigma2c} 
K\circ \Sigma_2^c = \Sigma_2^c \circ K\;,
\end{equation}
hence if $u$ is a K-real function, $\Sigma_2^c u$ is also a
K-real eigenfunction.\\
 We now show
\begin{proposition}\label{prop52}~\\ 
If $\dot \Omega$ has the symmetries of the
  rectangle \eqref{symrect}, then 
the multiplicity
 of the groundstate energy of $-\Delta_\Ab$ is $2$.\\
More generally the multiplicity of any eigenvalue is even.
\end{proposition}
\begin{proof}~\\ 
As observed in \cite{HHOO}, 
we can reduce the analysis of the Aharonov-Bohm Hamiltonian
  to the $K$-real space $L^2_K$
 where
$$
L^2_K(\dot{\Omega})=\{u\in L^2(\dot{\Omega}) \;,\; Ku =u\,\}\,.
$$
The scalar product on $L^2_K$, making of $L^2_K$ a real Hilbert space, is obtained by restricting
the scalar product on $L^2(\dot \Omega)$  to $L^2_K$  and it is immediate to verify that
$\langle u\,,\, v \rangle$ is indeed real for $u$ and $v$ in $L^2_K$.\\
Observing now that 
\begin{equation}
\Sigma_1^c \circ \Sigma_1^c =I\;,
\end{equation}
we obtain by writing
$$
u = \frac 12 (I+ \Sigma_1^c) u +  \frac 12 (I- \Sigma_1^c) u \,,
$$
an orthogonal  decomposition of $ L^2_K$ into
\begin{equation}\label{decomp1}
L^2_K = L^2_{K,\Sigma_1}\oplus
L^2_{K,a\,\Sigma_1}\,,
\end{equation}
where 
$$
 L^2_{K,\Sigma_1} =\{u\in L^2_K\;,\; \Sigma_1^c u=u\;\},
$$
and
$$
 L^2_{K,a\,\Sigma_1} =\{u\in L^2_K\;,\; \Sigma_1^c u=-u\;\}.
$$
We have  just to show that the restriction $\Pi_1$ of  $\frac 12 (I+ \Sigma_1^c)$
 to $L^2_K$ 
\beq
\Pi_1 := \frac 12 (I+ \Sigma_1^c)_{/ L^2_K}\,,
\eeq
is a projector. It is indeed clear that $\Pi_1$ is ($\mathbb R$-)linear 
and that $\Pi_1^2 =\Pi_1\,.$
It remains to verify that $\Pi_1^*=\Pi_1$. But we have, for $u$, $v$ in $L^2_K$,
$$
\langle \Sigma_1^c u\,,\, v\rangle
= i\langle \Gamma v \,,\, \Sigma_{1}u\rangle
= i\langle \Gamma\Sigma_{1} v\,,\, u\rangle
 =  \langle \Sigma_1^c v\,,\, u\rangle =  \langle u\,, \Sigma_1^c
 v\rangle\,.
$$
Moreover the  decomposition \eqref{decomp1}  is respected by $-\Delta_{\Ab}$.\\
Similarly, one can define the projector $\Pi_2$ by restriction of
 $\Sigma_2^c$ to $L^2_K$.\\
The second statement of Proposition \ref{prop52} will be a 
consequence of the following lemma
\begin{lemma}\label{lSigma3}~\\
Let 
\begin{equation}
\Sigma_3^c =\Sigma_1^c \;\Sigma_2^c\;,
\end{equation}
then $\Sigma_3^c$ commutes with  $-\Delta_{\Ab}$
 and $\Pi_3:=(\Sigma_3^c)_{/L^2_K}$ is a unitary operator from 
 $ L^2_{K,\Sigma_1}$
 onto  $ L^2_{K,a\,\Sigma_1}$.
\end{lemma}
\begin{proof} We note that 
\begin{equation}
\Sigma_3^c = i\, \Sigma_3\,,
\end{equation} 
 where $\Sigma_3$ is associated with:
\begin{equation}
\sigma_3(x,y) = (-x,-y)\,.
\end{equation}
The lemma follows then from the property that
 if $u$ is a solution of $Ku =u$ and $\Sigma_1^cu=u$, then
$$
\Sigma_1^c\; \Sigma_3^c \,u = \Sigma_1^c \; \Sigma_1^c \;\Sigma_2^c\, u
 = -  \Sigma_1^c \; \Sigma_2^c \; \Sigma_1^c\, u = - \Sigma_3^c\, u\;,
$$
where we have used the anticommutation of $\Sigma_1^c$ and
 $\Sigma_2^c$:
\begin{equation}\label{anticom}
\Sigma_1^c\; \Sigma_2^c = - \Sigma_2^c\;  \Sigma_1^c\,.
\end{equation}
\end{proof}
It remains to show that the first eigenvalue has multiplicity $2$. We
already know that 
it has an even multiplicity. It is enough to prove that  the
multiplicity  is at most $2$. Here 
we can 
 use   the results of \cite{HHOO}. Actually those results have to  be
 extended slightly\footnote{The paper of 
Alziary-Fleckinger-Takac \cite{AFT} is dealing with this case.} since  they are obtained assuming that the domain is homeomorphic to an 
annulus (so we are in a limiting case). It has been shown in
\cite{HHOO} that 
the nodal  set of a K-real 
groundstate is a line joining the center to the outer boundary. If the 
multiplicity of the groundstate eigenvalue is strictly greater than
2 we can, as in \cite{HHOO}, construct by linear combination of eigenfunctions
 a groundstate for which the zeroset hits the outer boundary at
two distinct points, hence a contradiction.
\end{proof}

We observe that the proof of the proposition gives more explicitly
the decomposition of $-\Delta_{\Ab}$ on $L^2_K$ 
into the direct orthogonal sum of
two unitary equivalent Hamiltonians. What we have done with
$\Sigma_1^c$ can be similarly done with $\Sigma_2^c$. This gives
immediately the following proposition.
\begin{proposition}\label{pisospectral}~\\
The four following operators $-\Delta_{\Ab, \Sigma_1}$,
  $-\Delta_{\Ab, a\,\Sigma_1}$,  $-\Delta_{\Ab, \Sigma_2}$
 and  $-\Delta_{\Ab, a\,\Sigma_2}$
 respectively defined  by the restriction of $-\Delta_{\Ab}$
to $ L^2_{K,\Sigma_1}$, $L^2_{K,a\,\Sigma_1}$, 
$L^2_{K,a\,\Sigma_2}$ and $L^2_{K,a\,\Sigma_2}$ 
are isospectral to $-\Delta_{\Ab}$. Moreover $\lambda$ is an
eigenvalue of any of the first four operators with multiplicity
$k(\lambda)$ if and only if $\lambda$ is an eigenvalue of 
multiplicity $2 k(\lambda)$ of $-\Delta_{\Ab}\,$.
\end{proposition}

Now we would like to analyze the nodal patterns 
 of eigenfunctions in the various symmetry spaces.
\begin{lemma}~\\
If $u\in C^\infty(\dot{\Omega}) \cap L^2_{K,\Sigma_2} $
  then its nodal set  contains  $ [-\frac a2, 0]\times \{0\}$. 
 Moreover, in $\dot{\Omega} \setminus \{ ]-\frac a2, 0]\times \{0\}\}\,$, 
  $$v=\exp- \frac{i\theta}{ 2}\; u$$ satisfies $$ \Gamma v=v\;\mbox{ and }
 \Sigma_2 \,v=v\,.$$ 
 \end{lemma}
\begin{proof}~\\
Noting that $Ku=u$ and $\Sigma_2^c u=u$ for $y=0$ and $x<0$ we
immediately obtain that $ \overline{u}(x,0)=u(x,0)=-
\overline{u}(x,0)$. Hence $u(x,0)=0$ for $x<0$.\\
An immediate computation gives 
$$
\Gamma v = \exp  \frac{i\theta}{2} \;\Gamma u = \exp -  \frac {i\theta} 2\;
K u = \exp -\frac {i\theta} 2\; u = v\,,
$$
and
$$
\Sigma_2 v =  \exp  \frac{i\theta}{2}\; \Sigma_2 u
 = \exp  \frac{i\theta}{2} \; \Gamma u  =v\,.
$$
\end{proof}
We would like to compare to some Dirichlet-Neumann problems 
on half-domains. We call upperhalf-$\Omega$ the set
\begin{equation}
\Omega^{uh}= \Omega \cap \{y>0\}\;,
\end{equation}
and introduce similarly the lowerhalf, lefthalf and righthalf domains
defined by
\begin{equation}
\Omega^{lh} =\Omega \cap \{y<0\}\;,\; 
\Omega^{leh}=\Omega \cap \{x<0\}\;,\;
\Omega^{rih}=\Omega \cap \{x>0\}\,.
\end{equation}
\begin{figure}[h!t]
\begin{center}
\hfill \includegraphics[height=3cm]{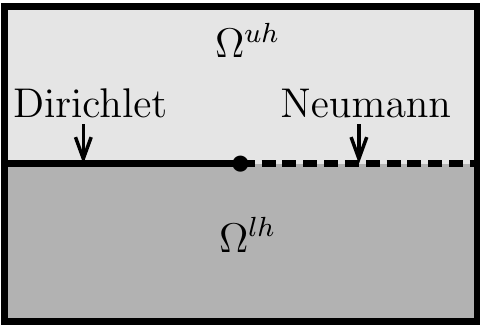} \hfill
\includegraphics[height=3cm]{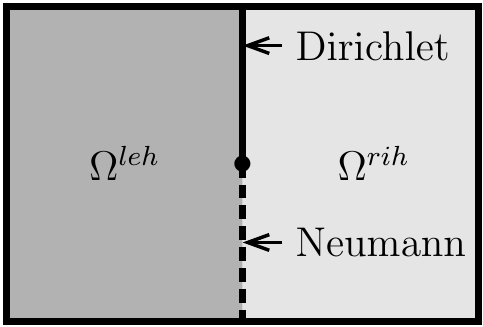} \hfill \ 
\caption{Domains $\Omega^{uh}$, $\Omega^{lh}$, $\Omega^{leh}$, $\Omega^{rih}$.}
\end{center}
\end{figure}
The previous lemma leads  to
\begin{proposition}~\\
If $u$ is a K-real $\Sigma_2^c$ invariant eigenfunction of
$-\Delta_\Ab$, then  the restriction to $\Omega^{uh}$
 of $\exp - \frac {i\theta}  2 \, u$
 is a real  eigenfunction of the   realization of  the Laplacian
 in $\Omega^{uh}$,
 with the following 
 Dirichlet-Neumann condition at $\pa \Omega^{uh}$: 
Dirichlet except on $ ]0,\frac a2]\times \{0\}$
 where we put Neumann.\\
In particular, if $\lambda$ is an eigenvalue of $-\Delta_{\Ab}$,
 then $\lambda$ is an eigenvalue of the Laplacian in $\Omega^{uh}$
 with  this Dirichlet-Neumann condition.
 \end{proposition}
Conversely, if we consider an eigenfunction $u$ of the Laplacian in $\Omega^{uh}$ 
 with this Dirichlet-Neumann condition
  and extend it into $ u^{ext}$  by symmetry in $\Omega \setminus ]-\frac
 a2,0]\times \{0\}$, then 
\begin{equation}\label{defv}
v = e^{\frac {i\theta} 2 } \; u^{ext}
\end{equation}
is a K-real eigenfunction of the Aharonov-Bohm Laplacian. 
More precisely, the function $v$ is first defined by the formula 
\eqref{defv}  for $\theta \in ]-\pi,\pi[$ and then extended as a 
$L^2$ function on the punctured square. Due to the properties of $u$, one can verify
 that $v$ is in the form domain of the Aharonov-Bohm operator.
Starting from a groundstate of the $DN$-problem in $\Omega^{uh}$, we get 
a $K$-real eigenstate in $L^2_{K,\Sigma_2}$ of the Aharonov-Bohm operator.\\

One can do the same type of argument in the three other 
 half-domains and obtain the following result.
\begin{proposition}~\\
The following  problems have the same eigenvalues:
\begin{itemize}
\item The Dirichlet problem for the Aharonov-Bohm 
 operator on $\dot{\Omega}$.
\item The Dirichlet-Neumann problem for the Laplacian 
on  $\Omega^{uh}$.
\item  The Dirichlet-Neumann problem  for the Laplacian 
on   $\Omega^{leh}$.
\item The Dirichlet-Neumann problem for the Laplacian 
on  $\Omega^{lh}$.
\item  The Dirichlet-Neumann problem  for the Laplacian 
on  $\Omega^{rih}$.
\end{itemize}
\end{proposition}
Of course this applies in particular to the case of the rectangle.

Let us go a little further by giving explicitly a unitary operator from
$L^2_{K,\Sigma_1}$ onto  $L^2_{K,\Sigma_2}$ proving the 
isospectrality. This is the object of 
\begin{lemma}~\\
The operator $$ U_{21}:=\frac{1}{\sqrt 2} (I+ \Sigma_2^c) $$
is a unitary operator from $L^2_{K,\Sigma_1}$ 
onto  $L^2_{K,\Sigma_2}$ whose inverse is given by
$$
 U_{12}:=\frac{1}{\sqrt 2} (I+ \Sigma_1^c)\,.
 $$
\end{lemma}
\begin{proof} 
Let $u\in L^2_{K,\Sigma_1}$, then, using \eqref{anticom}, 
$$
U_{12} U_{21} u = \frac 12 (I + \Sigma_1^c + \Sigma_2^c + \Sigma_1^c
\Sigma_2^c)u
=  \frac 12 (I + \Sigma_1^c + \Sigma_2^c - \Sigma_2^c
\Sigma_1^c)u =u\,.
$$
The proof that $ U_{21} U_{12} = I$ on $
L^2_{K,\Sigma_2}$
is obtained in the same way. Let us prove the norm conservation.
If $u\in L^2_{K,\Sigma_1}$, then
$$
|| U_{21} u||^2 = ||u||^2 + \frac 12 \langle \Sigma_2^c u\,,\,   u \rangle +
 \frac 12 \langle u\,,\, \Sigma_2^c u\rangle \,.
$$
But if $\Sigma_1^c u=u$, we can write
\begin{eqnarray*}
 \langle u\,,\,  \Sigma_2^c  u \rangle
& =&  \langle u\,,\, \Sigma_2^c \Sigma_1^c  u \rangle\\
&=& i \langle u\,,\,  \Sigma_2\Sigma_1 u \rangle\\
&=& i \langle \Sigma_2 \Sigma_{1}u\,,\,  u \rangle\\
&=& -\langle \Sigma_{2}^c u\,,\, \Sigma_{1}^c u \rangle\\
& = &- \langle \Sigma_2^c u\,,\, u\rangle \,.
\end{eqnarray*}
This leads to
\begin{equation}
|| U_{21} u||^2 = ||u||^2\,,\, \forall
 u \in  L^2_{K,\Sigma_1}\,.
\end{equation}
\end{proof}

\section{Application to minimal $3$-partitions\label{sec.isosquare}}

\subsection{Discussion on the square}\label{sssq}
We   look at the first excited  eigenvalue of the Dirichlet problem in
 the punctured square. The rules of \cite{HHOO}
 give  constraints about the nodal structure
 of the $K$-real eigenfunctions, which were already used. In particular we should have
 an odd number of lines arriving at the center.
So it is clear that $\{0\}$ belongs to the nodal set.
If three lines arrive at $0$ and if the nodal partition is a 
$3$-partition of type (a) it is reasonable to hope that 
this will give us a minimal $3$-partition.\\

\begin{figure}[h!t]
\begin{center}
\hfill\includegraphics[width=2cm]{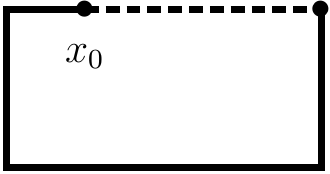}\hfill
\includegraphics[width=2cm]{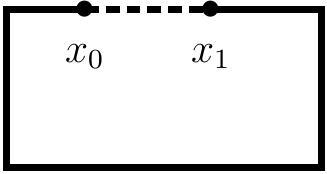}\hfill
\includegraphics[width=2cm]{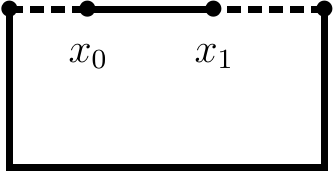}\hfill \ 
\caption{Mixed problems with Dirichlet-Neumann conditions.\label{fig.pbmixtes}}
\end{center}
\end{figure}

Let us explain the numerical strategy developed in \cite[\S 3]{BHV} to
exhibit a  
candidate for the minimal $3$-partition of the square. According to Theorem~\ref{partnod},
if the minimal $3$-partition of the square is admissible, it is associatd to the nodal set of an eigenfunction for $\lambda_{3}$. But there is no such function and therefore the minimal $3$-partition is non bipartite. 
Then we look for non bipartite  $3$-partitions whose topologies are enumerated in Proposition~\ref{mp3} and 
illustrated by Figure~\ref{fig.typeabc}. 
We first use the axial symmetry along the axis $\{y=0\}$. To recover a
partition of type (a), (b) or (c), we compute the second eigenfunction
and the next ones of  the mixed Dirichlet-Neumann Laplacian in $\Omega^{lh}$ with Dirichlet conditions except respectively on 
\begin{itemize}
\item $[x_{0},a/2]\times\{0\}$ for type (a);
\item $[x_{0},x_{1}]\times\{0\}$ for type (b);
\item $[-a/2,x_{0}]\times\{0\}\cup [x_{1},a/2]\times\{0\}$ for type (c).
\end{itemize}
These boundary conditions are illustrated in
Figure~\ref{fig.pbmixtes}. We move the points $x_{0}$ and $x_{1}$
along  the segment $[-a/2 ,a/2]\times\{0\}$. We expect to find an  eigenfunction such that, after symmetrization, their associated nodal sets constitute a $3$-partition and the nodal lines meet at the interior critical point with an angle of $2\pi/3$. \\

We recall the guess inherited  from numerical computations. When
minimizing over Dirichlet-Neumann problems in 
 $\Omega^{lh}$ (by putting Dirichlet except Neumann on
 $]x_{0},\frac a2]\times \{0\})$, one observes that  the minimal  second
     eigenvalue (as a function of $x_{0}$) such that  the two nodal domains 
  give rise by symmetry to a $3$-partition is obtained for  $x_{0}=0$.\\
Doing the same computation on the diagonal, we observe  also numerically
 that the  minimal second eigenvalue (as a function of the point on the 
 diagonal) such that the two nodal domains give rise by symmetry to 
 a three partition is obtained at the center.\\
Admitting these two numerical results, what we have proved for the
square is that the two minimal second eigenvalues are equal.\\
This could also suggest that we have a continuum 
of minimal $3$-partitions.\\
This point is not completely clear but could result of 
the analysis of the singularity at $\{0\}$.\\

When minimizing over Dirichlet-Neumann-Dirichlet or Neumann-Dirichlet-Neumann 
problems in $\Omega^{lh}$, the numerical computations (see \cite[\S 3.3]{BHV}) suggest 
that the nodal sets of the second eigenfunction never create a $2$-partition of $\Omega^{lh}$ 
leading by symmetry to a $3$-partition of $\Omega$. The corresponding 
 eigenmodes lead to too high energy,  hence do not qualify as possible
candidates for minimal 3-partitions.

\begin{remark}~\\
Note that we assumed  that the minimal 
partition is symmetric with respect to an axis of symmetry; the 
numerical experiments make this assumption plausible. 
One cannot a priori exclude
 that the first excited  $K$-real eigenfunction of the Aharonov-Bohm
 hamiltonian 
 consists  of one line joining $0$ to $\pa \Omega$
  and another line joining in $\dot \Omega $ 
 two points  of $\pa \Omega$. 
\end{remark}

\subsection{The symmetries of the square}

We now consider a convex domain which in addition to the invariance by
 $\sigma_1$ and $\sigma_2$ has
 an invariance by rotation (centered at the origin)  $r_{\frac \pi 2}$ of $\frac \pi 2$. We have typically in mind
  the case of the square. This
rotation  can be
quantized
 by
\begin{equation}\label{defrot}
\mathcal R_{\frac \pi 2} u (\cdot)= u ( r_{-\frac \pi 2} \,\cdot\,)\;,
\end{equation} 
where  $r_\alpha$ is the rotation by $\alpha$ in the plane.  
We  observe from \eqref{BAsecondform}  that this quantization commutes
with the operator:
\begin{equation}\label{commut2}
\Delta_{\Ab} \mathcal R_{\frac \pi 2} = \mathcal R_{\frac \pi 2} 
\Delta_{\Ab}\;,
\end{equation}
(and with its Dirichlet realization in  $\dot{\Omega}$). \\
More generally, we have the following Lemma
\begin{lemma}~\\
If $u$ is a K-real  eigenfunction of the Aharonov-Bohm  Hamiltonian,
 then  $u$ and $e^{i\frac \pi 4} \mathcal R_{\frac \pi 2}u$
 are linearly independent K-real eigenfunctions.
\end{lemma}
\begin{proof}
Let us first verify the K-reality.\\
We note that
$$
 \mathcal R_{\frac \pi 2}\, K  v = \mathcal R_{\frac \pi 2} \,\exp i
 \theta \, \Gamma v = \exp -i \frac \pi 2 \,\exp i
 \theta\, \Gamma \, \mathcal R_{\frac \pi 2}   v\;,
$$
hence
$$
 \mathcal R_{\frac \pi 2}\, K = \exp -i \frac \pi 2\,  K\,  \mathcal R_{\frac
 \pi 2}\,.
$$
This can be rewritten in the form
\begin{equation}
\left( \exp i\frac \pi 4 \, \mathcal R_{\frac \pi 2}\right) \circ  K
 = K \circ \left( \exp i\frac \pi 4 \,  \mathcal R_{\frac \pi
 2}\right)\,.
\end{equation}
This proves the first statement.\\
We now show that $\exp i\frac \pi 4 \, \mathcal R_{\frac \pi
  2} u$
 and $u$ are linearly independent (over $\mathbb R$) inside the real
  space 
$
L^2_K\,.
$
Let us look at the points of the nodal set belonging to
 the exterior boundary. Their cardinality should be odd by
 a result of \cite{HHOO} on K-real eigenfunctions.
If $u$ and $e^{i\frac \pi 4} \mathcal R_{\frac \pi 2}u$
 were proportional, this subset
  should be invariant by rotation of $\frac \pi 2$
 and should have consequently an even cardinality. Hence a contradiction.
\end{proof}
\begin{proposition}\label{pext}~\\
In the case of a convex domain having the symmetries of the square,  the Dirichlet-Neumann 
problem for the Laplacian  on the four  half-domains respectively defined by
\begin{equation}\begin{array}{ll}
\Omega^{--dh} =\Omega \cap \{ x+y <0\}\;,\; &
\Omega^{++dh} =\Omega \cap \{ x+y >0\}\;,\; \\
\Omega^{+-dh} =\Omega \cap \{ x-y >0\}\;,\; &
\Omega^{-+dh} =\Omega \cap \{ x-y <0\}\;,\;
\end{array}
\end{equation} 
are  also isospectral to the 
 problems introduced in Proposition \ref{pisospectral}.
\end{proposition}
\begin{figure}[h!t]
\begin{center}
\hfill \includegraphics[height=3cm]{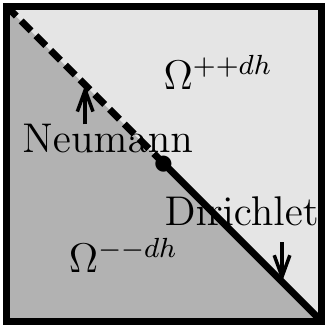} \hfill
\includegraphics[height=3cm]{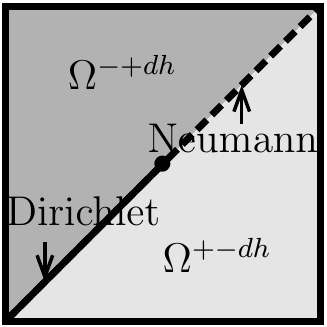} \hfill \ 
\caption{Domains $\Omega^{--dh}$, $\Omega^{++dh}$, $\Omega^{-+dh}$, $\Omega^{+-dh}$.}
\end{center}
\end{figure}

\begin{proof} We explain below  how to get the proposition  (which has also been
 verified numerically). 
We start from $v\in L^2_{K,\Sigma_2}$. 
Let us now consider
$$
w = v + \exp i\frac \pi  4\, \mathcal R_{\frac \pi 2} v\,.
$$
We have already shown that  $w$ is not zero and hence is  an 
 eigenfunction. It remains to analyze its zero set which 
  should contain  an half-diagonal.\\
Let us introduce 
$$
\Sigma_4^c = \exp i \frac \pi 4 \mathcal R_{\frac \pi 2} \Sigma_2^c\,.
$$
Using the property that
$$
\mathcal R_{\frac \pi 2} \Sigma_2^c = \Sigma_2^c  \mathcal R_{-\frac
  \pi 2}\;,
$$
we can verify that
\begin{equation}
\Sigma_4^c w = w\;,\; K w = w\,.
\end{equation}
For $\theta = - \frac{3 \pi}{ 4}$ and $x+iy =r \exp i \theta$, we get
$$
\overline{w} \exp i \frac \pi 4 = w\;,\; 
 \exp -i\frac { 3\pi}{ 4}  \overline{w} =w \;,
$$
hence $w=0$ for $\theta =- \frac{3\pi}{4}\;$. So the restriction
 of $w$ to $\Omega^{-+dh}$ multiplied by a phase factor leads to an
 eigenfunction of the DN-problem for the Laplacian in $\Omega^{-+dh}$.

The converse does not introduce new problems.
\end{proof}

\subsection{The covering approach}
As in \cite{HHOT}, we can also rewrite all the proofs by lifting the
  problem on the double covering 
  $\dot{\Omega}^\mathcal R$,   using the correspondence between
 $L^2_K$ and the real subspace of the functions in 
  $L^2(\dot{\Omega}^{\mathcal R})$ 
 such that $\Sigma u = -u$ where $\Sigma$ is associated with the map
  $\sigma$ by $(\Sigma u) (\omega)=u(\sigma(\omega))$. We recall with the
  notation of \cite{HH:2006}  that  $\sigma$
 is defined by associating  to each point $\omega $ of $\dot \Omega^{\mathcal R}$ the other point
 $\sigma(\omega)$ of 
$\dot \Omega^{\mathcal R}$ which has the same projection on $\dot
  \Omega $. Our initial proof was actually written in this way but we
  prefer to present in this paper another point of view.\\

In a recent paper, 
D. Jakobson, M. Levitin, N. Nadirashvili, I. Polterovich \cite{JLNP}
 obtain also nice isospectrality results involving Dirichlet-Neumann
 problems. They actually propose 
 three different proofs of their isospectrality results. The second one is a double covering argument which 
 is quite close to what we mentioned in the previous paragraph.
But there is no magnetic version and our magnetic examples seem to be new. So it would be
interesting  whether the magnetic approach can  produce
some isospectrality result which differ from the class of results
in this paper and the papers by Levitin, Parnovski and Polterovich
\cite{LPP} or \cite{PR}.

One should also mention that for these questions of isospectrality
 a  covering argument was already present in earlier works of
 B\'erard \cite{Be1,Be2,Be3}, B\'erard-Besson \cite{BeBe}, Sunada \cite{Su}....

\section{Some heuristics on the deformation of symmetric minimal
  partitions and numerical computations\label{sec.heuristics}}
\subsection{Some heuristics}
Here we discuss very heuristically in which general context the
numerical computations for the family of rectangles can be done.

One might investigate the special situation where tight upper and
lower bounds to  $\mathfrak L_3$ are available. We recall from
\eqref{compades3} that
$$ \la_k\le \mathfrak L_k\le L_k\,.$$
 If we have a family of domains depending analytically on a parameter 
$\alpha$, $\Om(\alpha)$, such that for some $k$,
\begin{equation}\label{lkLk}
\lim_{\alpha\downarrow 0}\la_k(\alpha)=L_k(0)\text{ and } \la_k(\alpha)<L_k(\alpha)\text{ for } \alpha>0\,,
\end{equation}
then we are led to the question how a minimal $k$-partition of
$\mathfrak D_k(\alpha)$ of $\Omega(\alpha)$ behaves as $\alpha$ tends to zero. In fact we might investigate 
more directly  $\mathfrak D_k(0)\,$. For $\alpha=0$ the only $\mathfrak
 D_k(0)$ are the one which are nodal partitions associated to $\la_k(0)=L_k(0)\,$.

We have  investigated this situation for rectangles. One has seen
that in this case the nodal partition can only have critical points at
the boundary.\\

Let us start with a slightly more general situation and consider a family $\Om(\beta)$ of simply connected domains 
which depends analytically upon a parameter $\beta\ge 0$. We assume that the spectrum and eigenfunctions corresponding to  $H(\Om(\beta))$ have for small $\beta\ge 0$ the following properties 
\begin{equation}\label{sigma0}
 \la_1(\beta)<\la_2(\beta)<\la_3(\beta)\le\la_4(\beta)=L_3(\beta)\,,
\end{equation}
with 
\begin{equation}\label{l3l4}
\la_3(\beta)<\la_4(\beta) \text{ for }\beta>0\text{ and }\la_3(0)=\la_4(0).
\end{equation}
Here as usual $H(\Om(\beta))$ is just $-\Delta$ with a Dirichlet boundary condition. 
We further assume that for $0<\beta<\beta_0$, $\mu(u_3(\beta))=2$, but that for $\beta=0$ 
the eigenspace $U_{3}$ of $\lambda_{1}(0)$ contains an eigenfunction $u\in U_3$ with $\mu(u)=3$. 
We hence have for $\beta>0$ that $\la_3(\beta)$ is not Courant-sharp and 
therefore there is a $\mathfrak L_3(\beta)>\la_3(\beta)$. 
But for $\beta=0$ we have $\mathfrak L_3(0)=\la_3(0)$.

\begin{figure}[h!t]
\begin{center}
\hfill \includegraphics[height=3cm]{FigBHHO08/figtypea1.pdf} \hfill
\includegraphics[height=3cm]{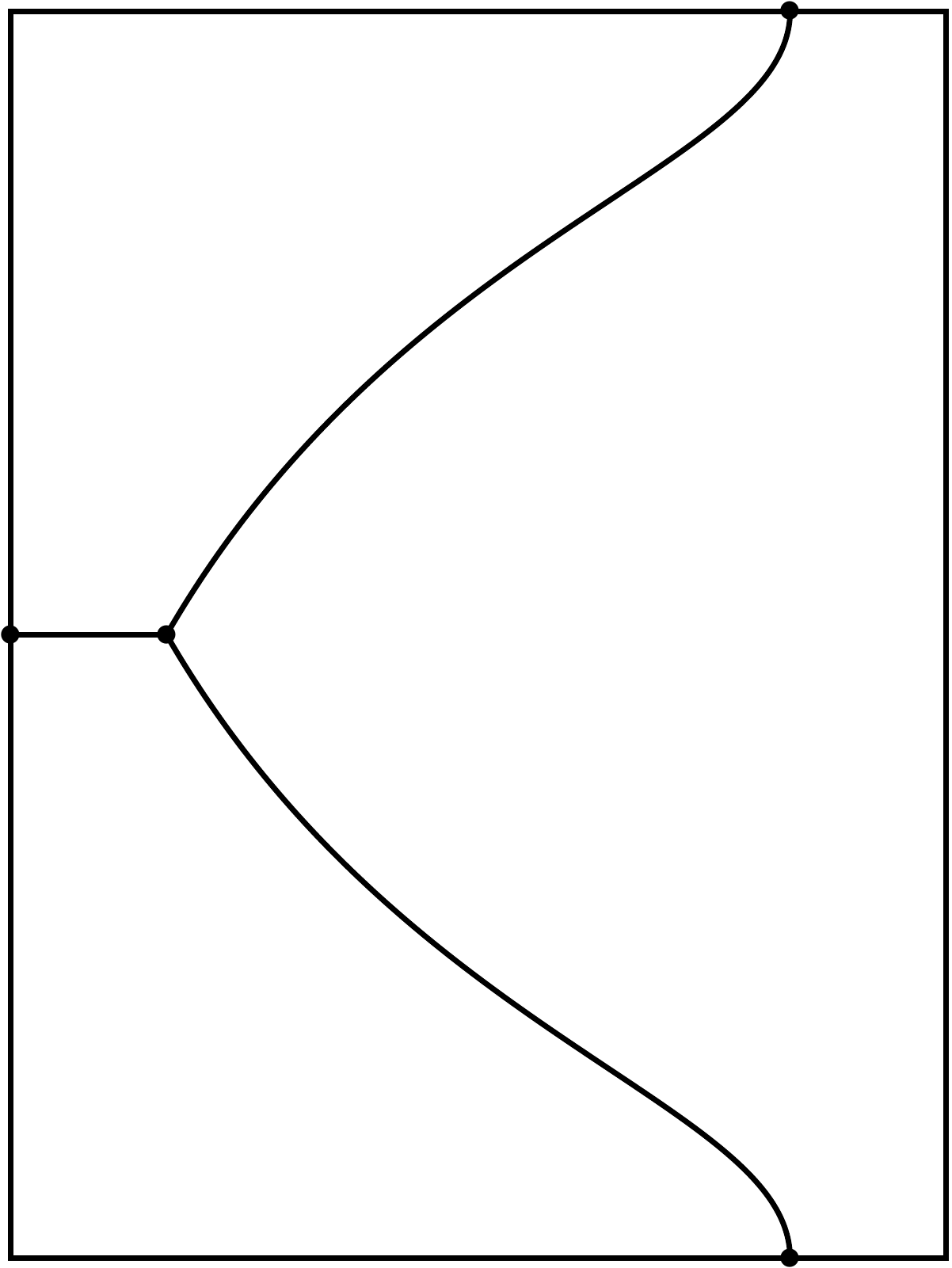} \hfill
\includegraphics[height=3cm]{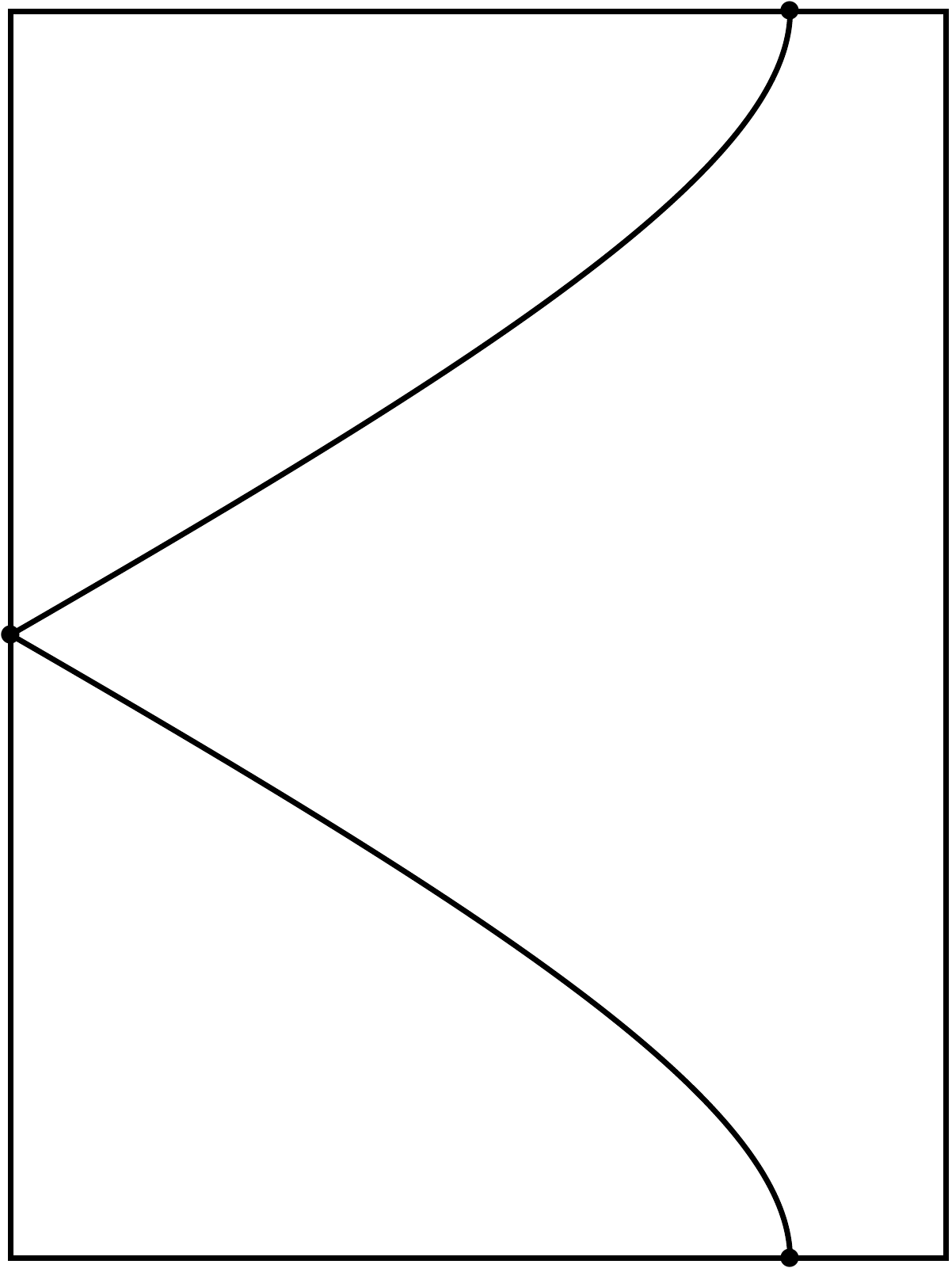} \hfill \ 
\caption{Deformation for type (a).\label{fig.typea}}
\end{center}
\end{figure}

\begin{figure}[h!t]
\begin{center}
\hfill \includegraphics[height=3cm]{FigBHHO08/figtypeb1.pdf} \hfill
\includegraphics[height=3cm]{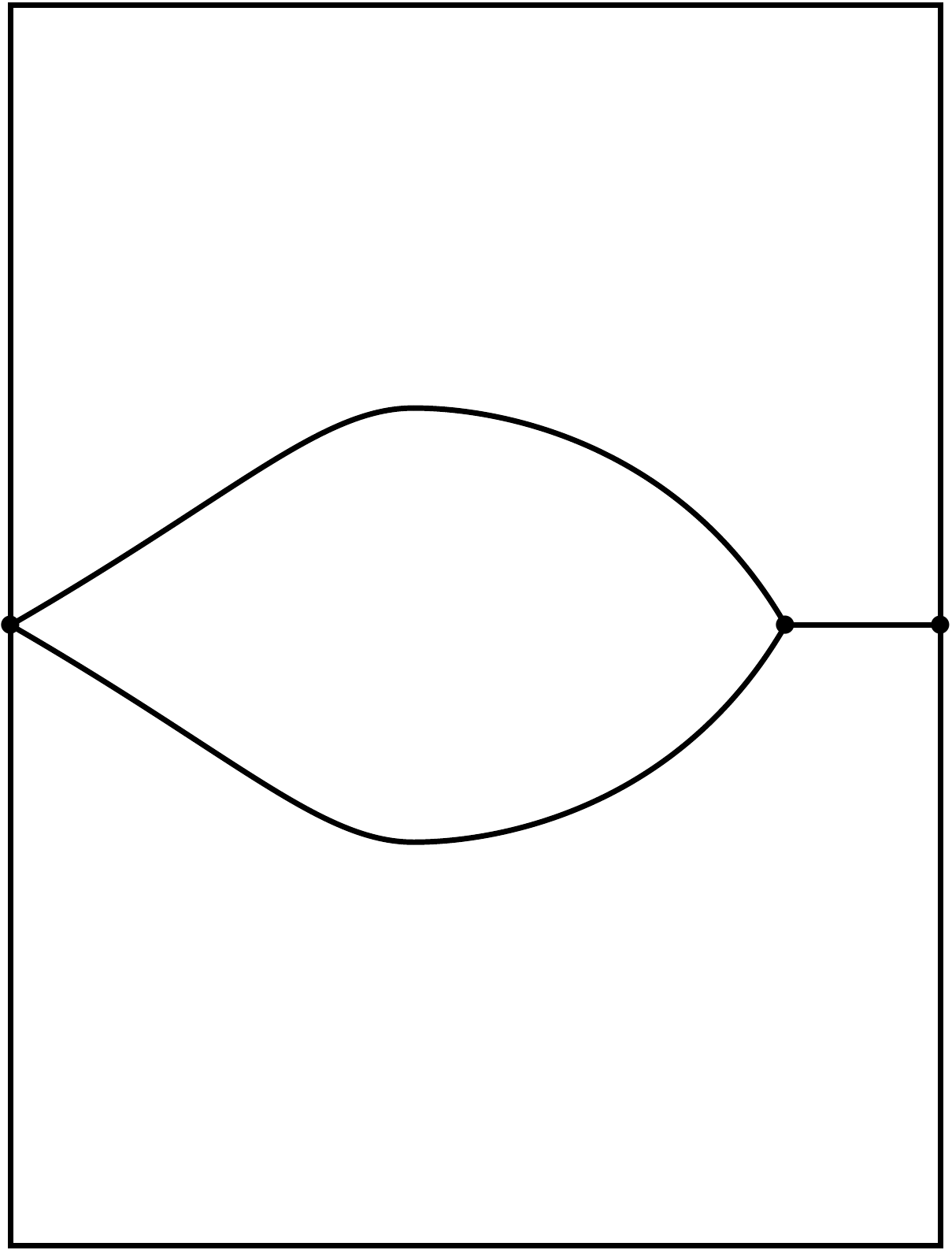} \hfill
\includegraphics[height=3cm]{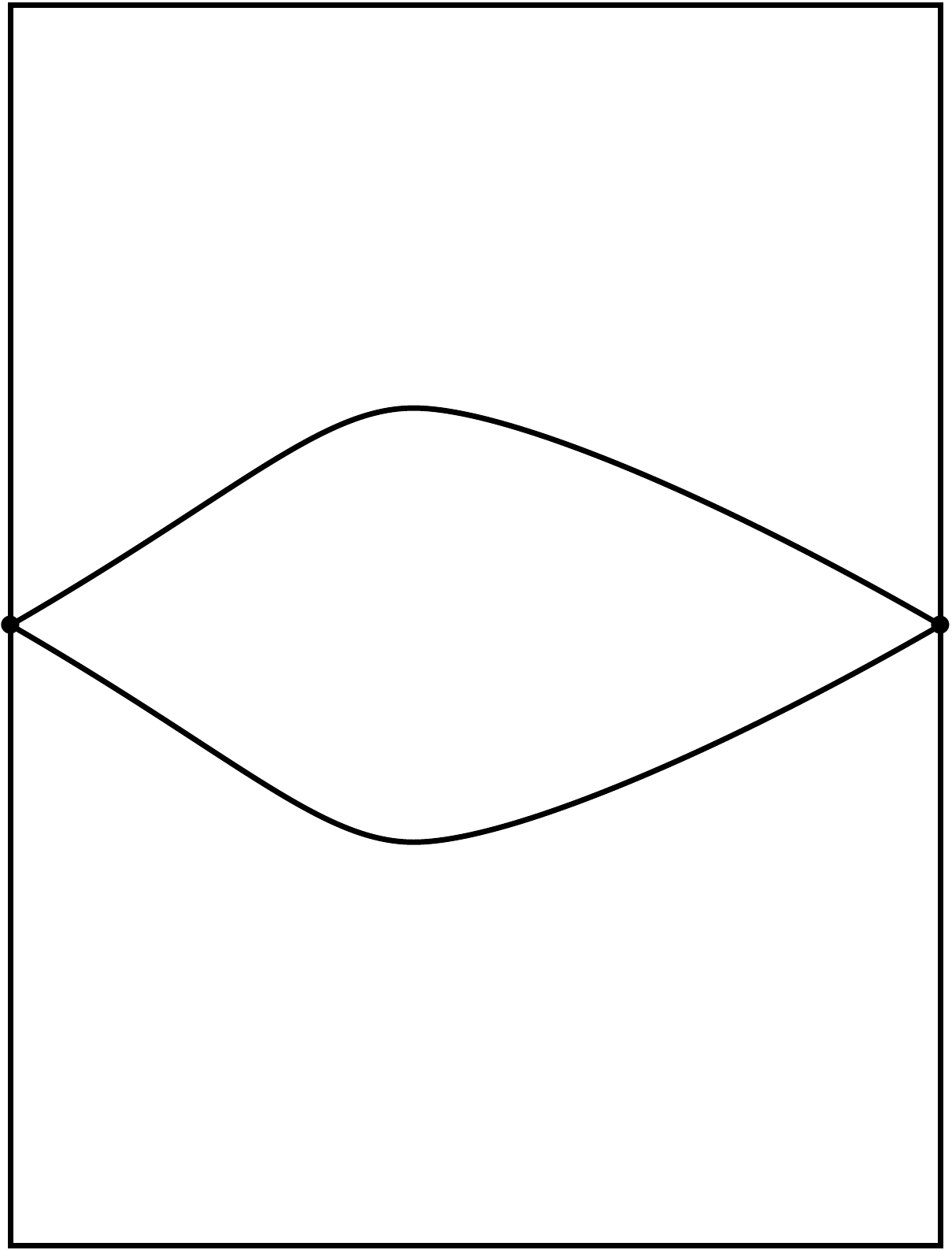} \hfill \ 
\caption{Deformation for type (b).\label{fig.typeb}}
\end{center}
\end{figure}

\begin{figure}[h!t]
\begin{center}
\hfill \includegraphics[height=3cm]{FigBHHO08/figtypec1.pdf} \hfill
\includegraphics[height=3cm]{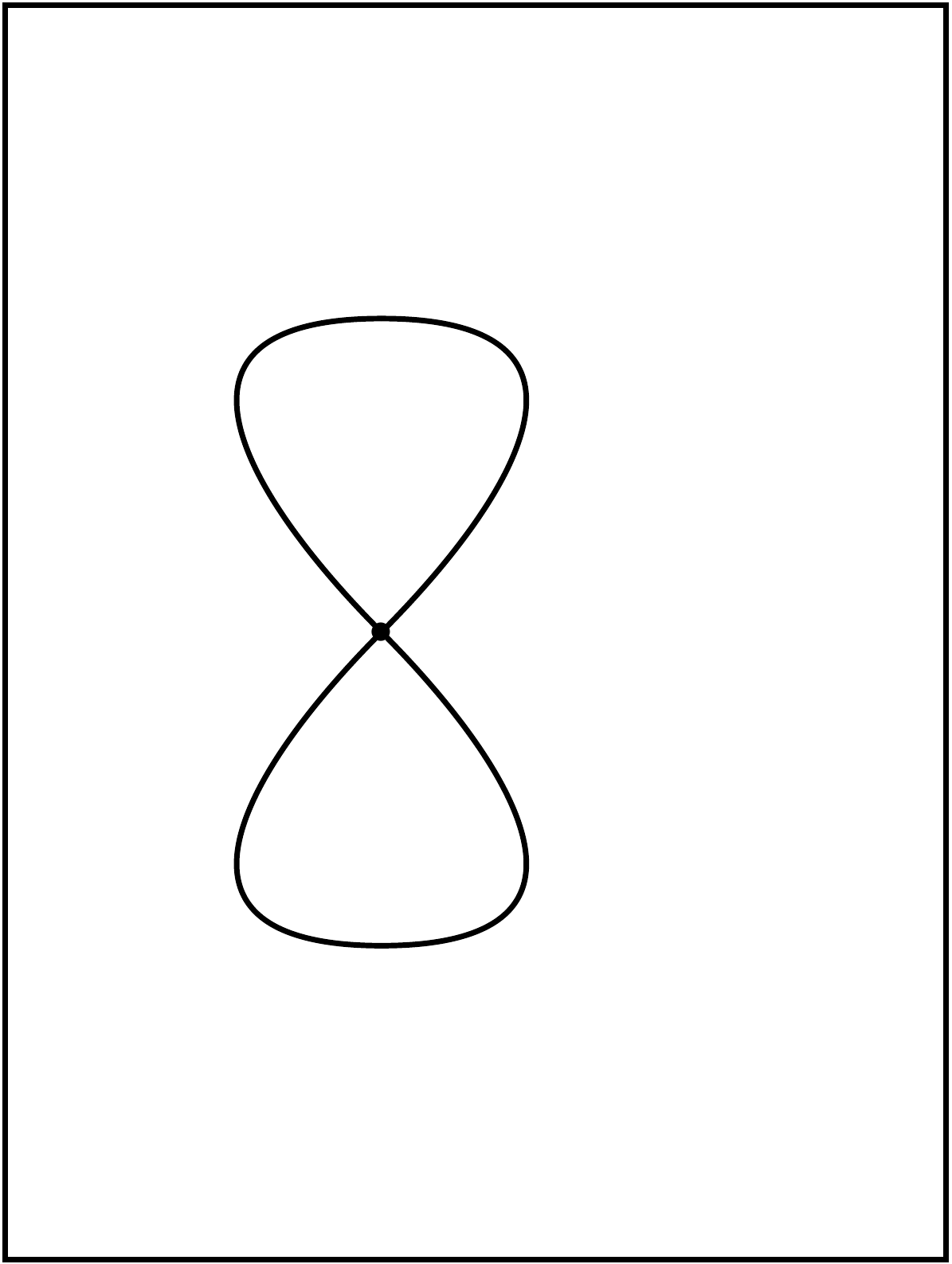} \hfill \ 
\caption{Deformation for type (c).\label{fig.typec}}
\end{center}
\end{figure}

According to Proposition \ref{mp3} we have for $\beta >0$ three types of
non bipartite partitions, \textbf{(a), (b), (c)}. 
But we have observed that for  $\mathfrak L_3(0)$, there is no
$\mathfrak D_3(0)$ which is not bipartite. So  we would like 
to understand how a non bipartite partition $\mathfrak D_3(\beta)$ can be 
deformed so that it becomes bipartite.  We emphasize that at the moment we have no
mathematical tools permitting to rigorously prove the validity of
these ``deformation argument'' but numerical tests show that they are
rather good for predicting what is observed.\\~\\
\noindent
\textbf{(a)} 
If the family is of type (a), what seems natural to imagine (see Figure~\ref{fig.typea}) 
is that the critical point should move to one point of the boundary and that 
(at least) two lines will start from this point and end at two other points of
 the boundary. These three points are not necessarily distincts.
\\\\\noindent 
\textbf{(b)}
If the family is of type (b), then the two critical points should again tend 
to the boundary (see Figure~\ref{fig.typeb}). We note indeed  that in the case when 
$z_1,z_2$ tend to a point  $z\in \Om$ then we get 4 nodal domains.     
\\\\\noindent
\textbf{(c)}
In the case (c), a new situtation can occur when the two critical points 
tend to one point in $\Omega$. One could indeed imagine a deformation of 
type (c) minimal partitions on a $3$-partition which is diffeomorphic to the 
figure eight (see Figure~\ref{fig.typec}).

Let us recall that in the case of the rectangle we have explicitly
verified that a limiting minimal partition cannot have a critical
point inside
 the rectangle.
\subsection{Numerics in the case for a family of rectangles}
\begin{figure}[h!t]
\begin{center}
\subfigure[$t=0$]{\includegraphics[height=4cm]{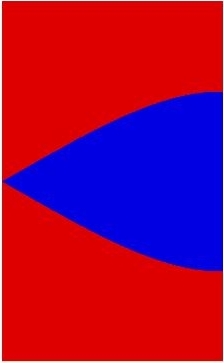} } \hfill
\subfigure[$t=2$]{\includegraphics[height=4cm]{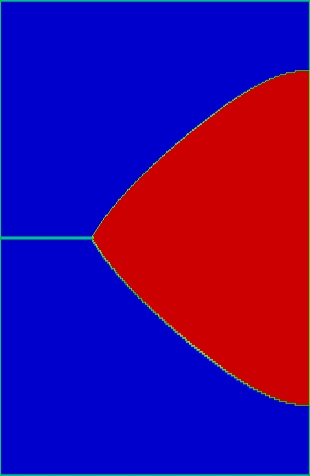} } \hfill
\subfigure[$t=3.45$]{\includegraphics[height=4cm]{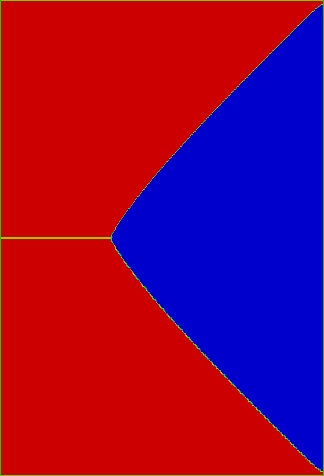} } \hfill
\subfigure[$t=7$]{\includegraphics[height=4cm]{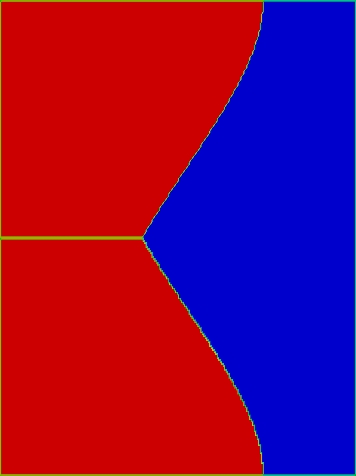}} \hfill
\subfigure[$t=11$]{\includegraphics[height=4cm]{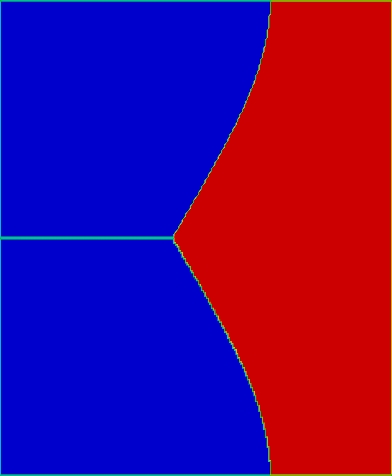}} \hfill
\subfigure[$t=15$]{\includegraphics[height=4cm]{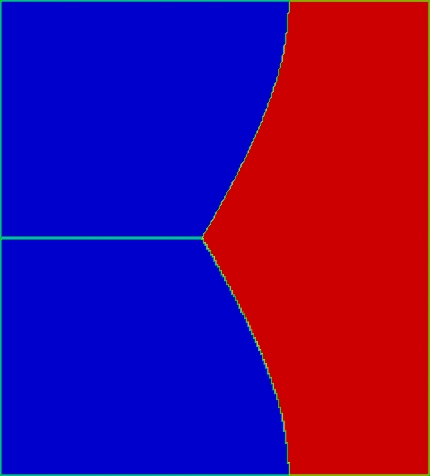}} \hfill
\subfigure[$t=20$]{\includegraphics[height=4cm]{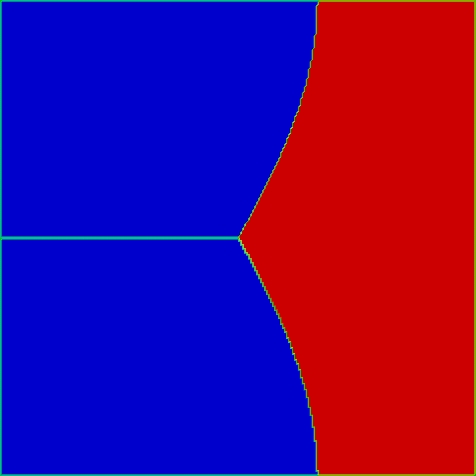}}
\caption{Simulations for rectangles $\Rb_{\pi\epsilon,\pi}$ with $\epsilon=(1-\frac{t}{20})\sqrt{\frac 38}+\frac{t}{20}$.\label{fig.recteps}}
\end{center} 
\end{figure}

We look at the case when the parameter $\beta$ introduced in the
 heuristic 
 subsection is the $\epsilon$
 of the computations and we go from $\epsilon =\sqrt{3/8}$  to
 $\epsilon =1$. We assume that all the minimal partitions are
 symmetric with respect to the horizontal axis $\{y=0\}$
 and of type (a). This permits to use the argument of reduction to an
 half-rectangle and to use the approach of the Dirichlet-Neumann
for each value of $\epsilon$. Figures~\ref{fig.recteps} present the evolution of the candidate to be minimal $3$-partition for rectangles $\Rb_{\pi\epsilon,\pi}$ with $\epsilon$ from $\sqrt{3/8}$ to $1$. \\
As mentioned in \cite{BHV}, numerical simulations on the half-square with mixed condition Dirichlet-Neumann-Dirichlet or Neumann-Dirichlet-Neumann never produce any $3$-partitions of type (b) or (c).  \\

\noindent{\bf Acknowledgements.}~\\
We thank M.~Van den Berg for indicating us the existence of
\cite{CyBaHo} and D.~Mangoubi for mentioning to one of us the recent
preprint 
\cite{PR}.
We thank also S.~Terracini for explaining us the results of 
her PHD Student B.~Norris on Aharonov-Bohm Hamiltonians
 and G.~Vial for discussion about the numerics. We also thank 
I.~Polterovich, D.~Jakobson, M.~Levitin, L.~Parnovski and 
G.~Verzini for useful discussions.


{\sc 
V. Bonnaillie-No\"el: IRMAR, ENS Cachan Bretagne, Univ. Rennes 1, CNRS, UEB,
av Robert Schuman, 35 170 Bruz, France.}\\
email: Virginie.Bonnaillie@bretagne.ens-cachan.fr

{\sc
B. Helffer: Laboratoire  de Math\'ematiques, Bat. 425,
Univ Paris-Sud and CNRS, 91 405 Orsay Cedex, France.}
\\
emai: Bernard.Helffer@math.u-psud.fr

{\sc 
T. Hoffmann-Ostenhof : International Erwin Schr\"odinger Institute for Mathematical Physics, Boltzmanngasse 9, A-090 Wien, Austria.} \\

email: thoffman@esi.ac.at

\end{document}